\documentclass[11pt]{amsart}

\usepackage{amsmath,amssymb}
\usepackage{amsfonts,amscd,amsxtra,calc,latexsym,mathtools}
\usepackage[mathscr]{eucal}
\usepackage{enumerate}
\usepackage{bm}
\usepackage[all]{xy}
\usepackage{tikz-cd}

\def\Spec{\mathrm{Spec}}
\def\Spf{\mathrm{Spf}}
\def\Hom{\mathrm{Hom}}
\def\Im{\mathrm{Im}}
\def\Ker{\mathrm{Ker}}
\def\Ext{\mathrm{Ext}}
\def\Coker{\mathrm{Coker}}
\def\mod{\mathrm{mod}}
\def\exp{\mathrm{exp}}
\def\Nil{\mathrm{Nil}}

\def\ZZ{\mathbb{Z}}
\def\QQ{\mathbb{Q}}

\def\FF{\mathbb{F}}
\def\GG{\mathbb{G}}

\def\XX{\mathbb{X}}
\def\YY{\mathbb{Y}}

\def\uu{\boldsymbol{u}}
\def\vv{\boldsymbol{v}}
\def\ww{\boldsymbol{w}}
\def\xx{\boldsymbol{x}}
\def\yy{\boldsymbol{y}}
\def\zz{\boldsymbol{z}}
\def\aa{\boldsymbol{a}}
\def\bb{\boldsymbol{b}}
\def\cc{\boldsymbol{c}}
\def\dd{\boldsymbol{d}}

\def\EE{\mathcal{E}}
\def\gg{\mathcal{G}}

\calclayout
\theoremstyle{plain}

\newtheorem{thm}{Theorem}[section]

\newtheorem{rem}[thm]{Remark}

\newtheorem{lem}[thm]{Lemma}
\newtheorem{sublem}[thm]{Sublemma}
\newtheorem{cor}[thm]{Corollary}

\numberwithin{equation}{section}
\numberwithin{table}{section}
\numberwithin{figure}{section}

\allowdisplaybreaks

\begin{document}

\title{Cartier Duals of Two-Dimensional Kummer--Artin--Schreier--Witt Type Kernels}
\author{\textbf{Michio Amano}}
\date{\today}
\address{School of Education, Meisei University, 2-1-1 Hodokubo Hino, Tokyo 191-8506, Japan}
\email{michio.amano@meisei-u.ac.jp}

\renewcommand\baselinestretch{1.1}

\begin{abstract}
We study Cartier duality for the kernel of a two-dimensional Kummer--Artin--Schreier--Witt type homomorphism over a $\mathbb{Z}_{(p)}$-algebra. For the $l$-th such homomorphism between extensions of deformations of the additive group by themselves, we give, under a natural nilpotency condition, an explicit description of the Cartier dual of its kernel in terms of Witt vectors. More precisely, the Cartier dual is identified with the kernel of a homomorphism between quotient fppf sheaves of two-dimensional Witt vectors determined by the operators naturally associated with the Kummer--Artin--Schreier--Witt type homomorphism. To prove this result, we use deformed Artin--Hasse exponential series and the formal completion of the isogeny.
\end{abstract}

\maketitle

\section{Introduction}

Let $A$ be a discrete valuation ring with maximal ideal $\mathfrak{m}$, and let $\lambda$ be a nonzero element of $\mathfrak{m}$. We define the $A$-group scheme $\mathcal{G}^{(\lambda)}_A=\Spec\, A[T,(1+\lambda T)^{-1}]$ whose group law is given by $x\cdot y = x+y+\lambda xy$. By Waterhouse and Weisfeiler \cite[Theorem~2.5]{WW}, any deformation of the additive group scheme $\mathbb{G}_a$ to the multiplicative group scheme $\mathbb{G}_m$ is of the form $\mathcal{G}^{(\lambda)}_A$. More precisely, let $\mathcal G$ be a smooth affine $A$-group scheme with connected fibers such that
$\mathcal{G}_K\simeq\mathbb{G}_{m,K}$ and $\mathcal{G}_k\simeq\mathbb{G}_{a,k}$, where $K=\mathrm{Frac}(A)$ and $k=A/\mathfrak{m}$. Then there exists a nonzero $\lambda\in\mathfrak{m}$, unique up to multiplication by a unit of $A$, such that $\mathcal{G} \simeq \mathcal{G}^{(\lambda)}_A$.

Let $p$ be a prime number, and let $\mathbb{Z}_{(p)}$ be the localization of the ring of integers $\mathbb{Z}$ at the prime ideal $(p)$. Let $\zeta$ be a primitive $p$-th root of unity, and set $\lambda:=\zeta-1$. Then $A:=\mathbb{Z}_{(p)}[\zeta]$ is a discrete valuation ring with uniformizer $\lambda$. Define a homomorphism of $A$-group schemes
\begin{align*}
\psi:\mathcal{G}^{(\lambda)}_A\rightarrow \mathcal{G}^{(\lambda^p)}_A;\ 
x \mapsto \psi(x):=\dfrac{1}{\lambda^{p}}\left\{(1+\lambda x)^p-1\right\}.
\end{align*}
Waterhouse \cite[Theorem~1]{W} and Sekiguchi, Oort and Suwa \cite[Section~2]{SOS} independently obtained the following Kummer--Artin--Schreier exact sequence:
\begin{equation*}
\begin{tikzcd}
0 \arrow[r] &
\left(\mathbb{Z}/p\mathbb{Z}\right)_A \arrow[r] &
\mathcal{G}^{(\lambda)}_A \arrow[r,"\psi"] &
\mathcal{G}^{(\lambda^{p})}_A \arrow[r] &
0.
\end{tikzcd}
\end{equation*}

Subsequently, Sekiguchi and Suwa described the extensions of $\mathcal{G}^{(\lambda_1)}_A$ by $\mathcal{G}^{(\lambda_2)}_A$ explicitly in terms of Witt vectors. Let $A$ be a discrete valuation ring of mixed characteristic $(0,p)$ with maximal ideal $\mathfrak m$, and let $\lambda_1,\lambda_2\in \mathfrak m\setminus\{0\}$. Let $W_A$ be the Witt group scheme over $A$, and let $F$ denote the Frobenius endomorphism of $W_A$. Put $F^{(\lambda)}:=F-[\lambda^{p-1}]$, where $[\lambda]$ denotes the Teichm\"uller lift of $\lambda\in A$. For a group scheme $G$, we write $\widehat{G}$ for the formal completion of $G$ along the zero section. Set
\begin{align*}
\widehat{W}(A_{\lambda_2})^{F^{(\lambda_1)}}
:=\Ker\!\left[F^{(\lambda_1)}:\widehat{W}(A_{\lambda_2})\to \widehat{W}(A_{\lambda_2})\right]
\end{align*}
where $A_{\lambda_2}:=A/(\lambda_2)$. Let $\iota:\Spec\,A_{\lambda_2} \hookrightarrow \Spec\, A$ denote the closed immersion. By \cite[Corollary~3.6]{S1}, there exists an isomorphism
\begin{align*}
\Ext^1\!\left(\mathcal{G}^{(\lambda_1)}_A,\mathcal{G}^{(\lambda_2)}_A\right)
\simeq
\Hom\!\left(\mathcal{G}^{(\lambda_1)}_A,\iota_\ast\mathbb{G}_{m,A_{\lambda_2}}\right)/\left\{(1+\lambda_1 x)^n \mid n\in\mathbb{Z}\right\}.
\end{align*}
Moreover, the isomorphism in \cite[Theorem~2.19.1]{SS1} is given by
\begin{align*}
\widehat{W}(A_{\lambda_2})^{F^{(\lambda_1)}}\xrightarrow{\sim}\Hom\!\left(\mathcal{G}^{(\lambda_1)}_A,\iota_\ast\mathbb{G}_{m,A_{\lambda_2}}\right);\ \uu\mapsto E_p(\uu,\lambda_1;X),
\end{align*}
where $E_p(\uu,\lambda_1;X)$ is the deformation of Artin--Hasse exponential series associated with $\uu$ (for details, see Subsection 3.2). Thus, the extension class $\left[\mathcal{E}_A^{(\lambda_1,\lambda_2;D)}\right]\in \Ext^1\!\left(\mathcal{G}^{(\lambda_1)}_A,\mathcal{G}^{(\lambda_2)}_A\right)$ is determined by a Witt vector $\uu\in\widehat{W}(A_{\lambda_2})^{F^{(\lambda_1)}}$. Here $D=D(X)\in A[X]$ is a normalized lift of $E_p(\uu,\lambda_1;X)\in A_{\lambda_2}[X]$, that is, a lift satisfying $D(0)=1$. Then the group scheme $\mathcal{E}^{(\lambda_1,\lambda_2;D)}_A$ is defined by
\begin{align*}
\EE^{(\lambda_1,\lambda_2;D)}_A=\Spec\, A\!\left[X,Y,\dfrac{1}{1+\lambda_1 X},\dfrac{1}{D(X)+\lambda_2 Y}\right]
\end{align*}
endowed with a group law for which the morphism
\begin{align*}
\alpha^{(\lambda_1,\lambda_2)}:\EE^{(\lambda_1,\lambda_2;D)}_A\to\mathbb{G}_{m,A}\times_{\Spec\,A}\mathbb{G}_{m,A};\ 
(x,y)\mapsto (1+\lambda_1 x,D(x)+\lambda_2 y)
\end{align*}
is a homomorphism of group schemes.

\begin{rem}\label{st:orig-1-1}
{\rm
Sekiguchi and Suwa constructed $\EE^{(\lambda_1,\lambda_2;D)}_A$ mainly over a discrete valuation ring. Next, M\'{e}zard, Romagny, and Tossici~\cite[Section~5]{MRT2} constructed $\EE^{(\lambda_1,\lambda_2;D)}$ over a $\ZZ_{(p)}$-algebra under the assumption that $\lambda_1$ and $\lambda_2$ are nonzero divisors.
}
\end{rem}

Let $l$ be a positive integer, and let $A$ be a $\mathbb{Z}_{(p)}$-algebra. Choose nonzero divisors $\lambda_1,\lambda_2\in A$ such that $\lambda_1^{p^l-p^{l-1}}\mid p$ and $\lambda_2^{p^l-p^{l-1}}\mid p$. Then, for each $0 \leq k \leq l-1$, we can take elements $\nu_{1,k},\nu_{2,k}\in A$ such that
\begin{align}\label{eq:orig-1-1}
\begin{split}
p^{l-k}=\nu_{1,k}\lambda_1^{p^l-p^k}\quad
\mbox{and}\quad
p^{l-k}=\nu_{2,k}\lambda_2^{p^l-p^k}
\end{split}
\end{align}
(for details, see Remark~\ref{st:orig-4-7}). Put
\begin{align*}
\boldsymbol{a}_1
:=\frac{1}{\lambda_1^{p^l}}\,p^l[\lambda_1]
=\left(\frac{b_0}{\lambda_1^{p^l}},\frac{b_1}{\lambda_1^{p^l}},\ldots\right)\in W(A),
\end{align*}
where $p^l[\lambda_1]=(b_0,b_1,\ldots)$. Note that, by \eqref{eq:orig-1-1}, $b_k/\lambda_1^{p^l}\in A$ for any $k\geq0$ (Remark~\ref{st:orig-4-9}). Take $\uu\in\widehat{W}(A_{\lambda_2})^{F^{(\lambda_1)}}$ and $\vv\in\widehat{W}(A_{\lambda_2^{p^l}})^{F^{(\lambda_1^{p^l})}}$. Then we can take a lift $\widetilde{\uu}\in\widehat{W}(A_{\lambda_2^{p^l}})$ of $\uu$ (Sublemma~\ref{st:orig-4-10}). Assume the congruence
\begin{align}\label{eq:orig-1-2}
p^l\widetilde{\uu}-j[\lambda_1]\equiv T_{\aa_1}(\vv)\pmod{\lambda_2^{p^l}},
\end{align}
where $T_{\aa_1}$ is a homomorphism of Witt group schemes defined in Subsection~2.2. By \eqref{eq:orig-1-1}, we can set the polynomial
\begin{align}\label{eq:orig-1-3}
\psi_1(X):=\dfrac{1}{\lambda_1^{p^l}}\left\{\left(1+\lambda_1 X\right)^{p^l}-1\right\}\in A[X].
\end{align}
Assume that $\lambda_1$ is nilpotent in $A_{\lambda_2}$. Let $D(X)\in A[X]$ be a normalized lift of $E_p(\widetilde{\uu},\lambda_1;X)\in A_{\lambda_2^{p^l}}[X]$, and let $D'(X)\in A[X]$ be a normalized lift of $E_p(\vv,\lambda_1^{p^l};X)\in A_{\lambda_2^{p^l}}[X]$. Thus $D(0)=D'(0)=1$. Then, by \eqref{eq:orig-1-2}, the congruence
\begin{align}\label{eq:orig-1-4}
D(X)^{p^l}(1+\lambda_1 X)^{-j}\equiv D'\left(\psi_1(X)\right)\pmod{\lambda_2^{p^l}}
\end{align}
holds (Lemma~\ref{st:orig-4-11}). By \eqref{eq:orig-1-1} and \eqref{eq:orig-1-4}, we can define the regular function
\begin{align}\label{eq:orig-1-5}
\psi_2(X,Y):=\dfrac{1}{\lambda_2^{p^l}}\left\{\dfrac{\left(D(X)+\lambda_2 Y\right)^{p^l}}{(1+\lambda_1 X)^j}-D'(\psi_1(X))\right\}\in A\!\left[X,Y,\dfrac{1}{1+\lambda_1X}\right].
\end{align}

\begin{thm}\label{st:orig-1-2}
Let $A$ be a $\mathbb{Z}_{(p)}$-algebra, and let $\lambda_1$ and $\lambda_2$ be nonzero divisors satisfying \eqref{eq:orig-1-1}. Assume that $\lambda_1$ is nilpotent modulo $\lambda_2$. Take $j\in\ZZ$. Choose $\uu\in \widehat{W}(A_{\lambda_2})^{F^{(\lambda_1)}}$ and $\vv\in \widehat{W}(A_{\lambda_2^{p^l}})^{F^{(\lambda_1^{p^l})}}$ satisfying \eqref{eq:orig-1-2}. Let $\mathcal{E}_A^{(\lambda_1,\lambda_2;D)}$ be the extension associated with $\boldsymbol{u}$, and let $\mathcal{E}_A^{(\lambda_1^{p^l},\lambda_2^{p^l};D')}$ be the extension associated with $\boldsymbol{v}$. Then the morphism defined by \eqref{eq:orig-1-3} and \eqref{eq:orig-1-5}
\begin{align*}
\Psi_{\uu,\vv}^{(l,j)}:\EE_A^{(\lambda_1,\lambda_2;D)}\to\EE_A^{(\lambda_1^{p^l},\lambda_2^{p^l};D')};\ 
(x,y)\mapsto \Psi_{\uu,\vv}^{(l,j)}(x,y):=(\psi_1(x),\psi_2(x,y))
\end{align*}
is a surjective $A$-group scheme homomorphism such that the diagram
\begin{equation}\label{eq:orig-1-6}
\begin{tikzcd}[column sep=12ex]
\mathcal{E}_A^{(\lambda_1,\lambda_2;D)} \arrow[r,"\alpha^{(\lambda_1,\lambda_2)}"]
\arrow[d, "\Psi^{(l,j)}_{\boldsymbol{u},\boldsymbol{v}}"'] &
\mathbb{G}_{m,A}^2 \arrow[d,"\Theta^{(l,j)}"]\\
\mathcal{E}_A^{(\lambda_1^{p^l},\lambda_2^{p^l};D')} \arrow[r,"\alpha^{(\lambda_1^{p^l},\lambda_2^{p^l})}"] &
\mathbb{G}_{m,A}^2
\end{tikzcd}
\end{equation}
commutes, where $\alpha^{(\lambda_1,\lambda_2)}(X,Y)=\bigl(1+\lambda_1 X,\; D(X)+\lambda_2 Y\bigr)$ and
\begin{align*}
\Theta^{(l,j)}:\mathbb{G}_{m,A}^2\rightarrow \mathbb{G}_{m,A}^2;\ 
(t_1,t_2)\mapsto\bigl(t_1^{p^l},\, t_1^{-j}t_2^{p^l}\bigr),
\end{align*}
with $\mathbb{G}_{m,A}^2:=\mathbb{G}_{m,A}\times_{\Spec A}\mathbb{G}_{m,A}$.
\end{thm}

\begin{rem}\label{st:orig-1-3}
{\rm
Take $l=j=1$. The homomorphism $\Psi_{\uu,\vv}^{(1,1)}$ defined in Theorem~\ref{st:orig-1-2} includes, as a special case, the homomorphism that gives rise to the Kummer--Artin--Schreier--Witt theory of degree $p^2$ (Remark~\ref{st:orig-4-12}).
}
\end{rem}

\begin{rem}\label{st:orig-1-4}
{\rm
The descriptions of $D(X)$ and $D'(X)$ in the case $l=1$ were studied in detail by Tossici \cite{To}. For filtered group schemes associated with models of $\mu_{p^n}$, see M\'{e}zard, Romagny, and Tossici \cite{MRT1}.
}
\end{rem}

Here we set
\begin{align*}
\mathcal{N}_A
:=\Ker\!\left[\Psi^{(l,j)}_{\boldsymbol{u},\boldsymbol{v}}:\mathcal{E}_A^{(\lambda_1,\lambda_2;D)}\rightarrow \mathcal{E}_A^{(\lambda_1^{p^l},\lambda_2^{p^l};D')}\right].
\end{align*}
Then $\mathcal{N}_A$ is a finite flat group scheme of order $p^{2l}$ over $A$. Our purpose is to describe the Cartier dual of $\mathcal{N}_A$ under a suitable restriction.
Let $\uu^{\#}\in W(A)$ be a lift of the chosen $\widetilde{\uu}\in\widehat{W}(A_{\lambda_2^{p^l}})$, and let $\widetilde{\vv}\in W(A)$ be a lift of $\vv$. We put
\begin{align*}
\bb:=\dfrac{1}{\lambda_2}F^{(\lambda_1)}(\uu^{\#})\in W(A)
\end{align*}
and we set
\begin{align*}
U_2:=
\begin{pmatrix}
F^{(\lambda_1)} & -T_{\boldsymbol{b}}\\
0 & F^{(\lambda_2)}
\end{pmatrix},
\end{align*}
which defines an endomorphism of $W_A^2:=W_A\times_{\Spec A} W_A$. Similarly, we set
\begin{align*}
\boldsymbol{b}'=\dfrac{1}{\lambda_2^{p^l}}F^{(\lambda_1^{p^l})}(\widetilde{\boldsymbol{v}})\in W(A)
\quad \text{and}\quad 
U'_2:=
\begin{pmatrix}
F^{(\lambda_1^{p^l})} & -T_{\boldsymbol{b}'}\\
0 & F^{(\lambda_2^{p^l})}
\end{pmatrix}.
\end{align*}
We put
\begin{align*}
\cc:=\dfrac{1}{\lambda_2^{p^l}}\left\{p^l\uu^{\#}-j[\lambda_1]-T_{\aa_1}(\widetilde{\vv})\right\}\in W(A).
\end{align*}
Moreover, we set
\begin{align*}
\boldsymbol{a}_2:=\dfrac{1}{\lambda_2^{p^l}}p^l[\lambda_2]\in W(A).
\end{align*}
We define the homomorphism
\begin{align*}
\Upsilon_2:=
\begin{pmatrix}
T_{\boldsymbol{a}_1} & T_{\cc}\\
0 & T_{\boldsymbol{a}_2}
\end{pmatrix}
\end{align*}
on $W^2_A$. Here we assume that $B$ is an $A$-algebra such that each $\nu_{i,k}$ is nilpotent. Let $C$ be a $B$-algebra. Then there exists a homomorphism $\Upsilon_2'$ such that the diagram
\begin{equation}\label{eq:orig-1-7}
\begin{tikzcd}
W_B^2(C) \arrow[r,"\Upsilon_2"] \arrow[d,"U_2'"'] &
W_B^2(C) \arrow[d,"U_2"]\\
W_B^2(C) \arrow[r,"\Upsilon'_2"] &
W_B^2(C)
\end{tikzcd}
\end{equation}
commutes (Lemma~\ref{st:orig-5-7}). We put
\begin{align}\label{eq:orig-1-8}
\left(W_B^2/\Upsilon_2\right)\!(C)
:=\Coker\!\left[\Upsilon_2:W_B^2(C)\to W_B^2(C)\right].
\end{align}
Moreover, we set
\begin{align}\label{eq:orig-1-9}
\left(W_B^2/\Upsilon'_2\right)\!(C)
:=\Coker\!\left[\Upsilon'_2:W_B^2(C)\to W_B^2(C)\right].
\end{align}
Note that $W_B^2/\Upsilon_2$ and $W_B^2/\Upsilon'_2$ are sheaves for the fppf topology on the category of $B$-algebras (Lemma~\ref{st:orig-5-10}). We consider the following diagram:
\begin{equation*}
\begin{tikzcd}
W_B^2(C) \arrow[r] \arrow[d,"U_2"'] &
\left(W_B^2/\Upsilon_2\right)\!(C)
\arrow[d,"\overline{U_2}"] \\
W_B^2(C) \arrow[r] &
\left(W_B^2/\Upsilon'_2\right)\!(C),
\end{tikzcd}
\end{equation*}
where $\overline{U_2}$ is defined by $\overline{U_2}(\overline{\boldsymbol{x}}):=\overline{U_2(\boldsymbol{x})}$. Then, by \eqref{eq:orig-1-7}, $\overline{U_2}$ is well-defined, and the above diagram commutes. We set
\begin{align*}
\mathcal{M}_B(C)
:=\Ker\!\left[\overline{U_2}:
\left(W_B^2/\Upsilon_2\right)\!(C)
\to
\left(W_B^2/\Upsilon'_2\right)\!(C)\right].
\end{align*}
Note that $\mathcal{M}_B$ is a sheaf for the fppf topology on the category of $B$-algebras. We will show that, for an $A$-algebra $B$ satisfying a suitable nilpotency condition of $\nu_{i,k}$, the Cartier dual of $\mathcal{N}_B$ is isomorphic to $\mathcal{M}_B$. In the formal exponential descriptions below, $D$ denotes $E_p(\uu^\#,\lambda_1;X)$ after the change-of-lift identification explained in Section~5. The main tool in our proof is the following isomorphism of Sekiguchi and Suwa \cite[Theorem~5.1]{SS2}:
\begin{align*}
\Ker\!\left[U_2:W_A^2(B)\to W_A^2(B)\right]
&\xrightarrow{\sim}
\Hom\!\left(\widehat{\mathcal{E}}^{(\lambda_1,\lambda_2;D)}_B,\widehat{\mathbb{G}}_{m,B}\right);\\
\xx&\mapsto E_p\left(\xx,(\lambda_1,\lambda_2);X,Y\right),
\end{align*}
where
\begin{align*}
E_p\left(\xx,(\lambda_1,\lambda_2);X,Y\right)
:=E_p\left(\xx_1,\lambda_1;X\right)\cdot E_p\!\left(\xx_2,\lambda_2;\frac{Y}{D(X)}\right)
\end{align*}
for $\boldsymbol{x}=(\boldsymbol{x}_1,\boldsymbol{x}_2)$ (for details, see Subsection~\ref{subsec:two-dimensional-exponentials}). Our main result is as follows.

\begin{thm}\label{st:orig-1-5}
Under the assumptions and notation of Theorem~\ref{st:orig-1-2}, let $B$ be an $A$-algebra. Assume that, for each $i=1,2$, every $\nu_{i,k}$ is nilpotent in $B$. Then, for any $B$-algebra $C$, the correspondence
\begin{align}\label{eq:orig-1-10}
\begin{split}
\mathcal{M}_B(C)
&\xrightarrow{\sim}
\left(\mathcal{N}_B\right)^D\!(C) \simeq\Hom\left(\mathcal{N}_C,\mathbb{G}_{m,C}\right);\\
\overline{\boldsymbol{x}}&\mapsto E_p\!\left(\boldsymbol{x},(\lambda_1,\lambda_2);\overline{X},\overline{Y}\right)\Gamma_{\boldsymbol{w}}\!\left(\overline{X},\overline{Y}\right)
\end{split}
\end{align}
defines an isomorphism of fppf sheaves on the category of $B$-algebras:
\begin{align*}
\mathcal{M}_B\simeq\left(\mathcal{N}_B\right)^D,
\end{align*}
where $\left(\mathcal{N}_B\right)^D$ denotes the Cartier dual of $\mathcal{N}_B$, $\boldsymbol{w}\in W^2_B(C)$ is the unique vector satisfying $U_2(\boldsymbol{x})=\Upsilon_2'(\boldsymbol{w})$, and $\Gamma_{\boldsymbol{w}}$ is the normalized correction cochain defined in \eqref{eq:gamma-definition} and Lemma~\ref{st:orig-5-2}. Here $(\overline X,\overline Y)$ are the coordinates on the kernel after the change of variables at the beginning of Section~5.
The existence of $\boldsymbol{w}$ follows from the definition of $\mathcal M_B(C)$, and its uniqueness follows from Lemma~\ref{st:orig-5-7}.
\end{thm}

\begin{rem}\label{st:orig-1-6}
{\rm
For the Frobenius specialization of Corollary~\ref{st:orig-5-13},
Theorem~\ref{st:orig-1-5} recovers the corresponding case of \cite[Theorem~2]{AA}.
}
\end{rem}

For the proof, we also use two general results. In
Theorem~\ref{st:orig-5-9}, we prove the vanishing of the positive
fppf cohomology of $W_R^r$ on $\Spec R$ for every finite $r$.
In Theorem~\ref{thm:boundary-cocycle}, we describe the boundary
homomorphism by a Hochschild cocycle when the quotient has global
formal coordinates. This gives the Hochschild exact sequence used
in the proof of Theorem~\ref{st:orig-1-5}.

The contents of this paper are as follows. In Sections~2 and~3 we
recall Witt vectors and deformed Artin--Hasse exponential series.
In Section~4 we construct the two-dimensional Kummer--Artin--Schreier--Witt
type isogeny and consider its formal completion. In Section~5 we
prove the two general results stated above and use them to prove
Theorem~\ref{st:orig-1-5}.
Throughout this paper, we use the following notations:
\begin{align*}
            \GG_{a,A}:\ &\textrm{additive group scheme over $A$},\\
            \GG_{m,A}:\ &\textrm{multiplicative group scheme over $A$},\\
  \widehat{\GG}_{m,A}:\ &\textrm{multiplicative formal group scheme over $A$},\\
                W_{A}:\ &\textrm{group scheme of Witt vectors over $A$},\\
                    F:\ &\textrm{Frobenius endomorphism of $W_{A}$},\\
                    V:\ &\textrm{Verschiebung endomorphism of $W_{A}$},\\
            [\lambda]:\ &\textrm{Teichm\"{u}ller lift $(\lambda,0,0,\ldots)\in W(A)$ of $\lambda\in A$},\\
        F^{(\lambda)}:\ &=F-[\lambda^{p-1}],\\
              T_{\aa}:\ &\textrm{homomorphism defined by $\aa \in W(A)$\ (recalled in Subsection~2.2)},\\
            \uu^{(p)}:\ &=(u_0^p,u_1^p,\ldots)\ \text{for}\ \uu=(u_0,u_1,\ldots)\in W(A),\\
        \Ext^1(G,H)\ :\ & \mbox{isomorphism classes of extensions of $G$ by $H$},\\
         H^2_0(G,H)\ :\ & \mbox{Hochschild cohomology of $G$ with coefficients in $H$}.
\end{align*}
Throughout, the category of $A$-algebras is endowed with the fppf topology.

\section{Witt vectors}

In this short section, we recall the facts about Witt vectors needed in this paper. For details, see \cite[Chap.~V]{DG} or \cite[Chap.~III]{HZ}.

\subsection{Definition of Witt vectors}

Let $\mathbb{X}=(X_0,X_1,\ldots)$ be a sequence of variables. For each $n\geq 0$, we denote by $\Phi_n(\mathbb{X})=\Phi_n(X_0,X_1,\ldots,X_n)$ the Witt polynomial
\begin{align*}
\Phi_n(\mathbb{X})=X_0^{p^n}+pX_1^{p^{n-1}}+\cdots+p^nX_n
\end{align*}
in $\ZZ [ \mathbb{X} ] = \ZZ [ X_0 , X_1 , \ldots ]$. Let
\begin{align*}
W_{n,\ZZ} = \Spec\, \ZZ [X_0,X_1,\ldots,X_{n-1}]
\end{align*}
be a $n$-dimensional affine space over $\ZZ$. The phantom map $\Phi$ is defined by
\begin{align*}
\Phi:W_{n,\ZZ}\rightarrow\mathbb{A}^n_\ZZ;\ 
\xx \mapsto (\Phi_0(\xx),\Phi_1(\xx),\ldots,\Phi_{n-1}(\xx)),
\end{align*}
where $\mathbb{A}^n_\ZZ$ is the usual $n$-dimensional affine space over $\ZZ$. The scheme $\mathbb{A}^n_\ZZ$ has a natural ring scheme structure. It is known that $W_{n,\ZZ}$ has a unique commutative ring scheme structure over $\ZZ$ such that the phantom map $\Phi$ is a homomorphism of commutative ring schemes over $\ZZ$. Then $A$-valued points $W_{n,\ZZ}(A)$ are called Witt vectors of length $n$ over $A$.
The restriction homomorphism $R_{n-1}:W_{n,\ZZ}\rightarrow W_{n-1,\ZZ}$ is given by
\begin{align*}
\xx=(x_0,x_1,\ldots,x_{n-1})\mapsto R_{n-1}(\xx):=(x_0,x_1,\ldots,x_{n-2}).
\end{align*}
Then the Witt scheme $W_{\ZZ}$ is defined by
\begin{align*}
W_\ZZ:=\varprojlim_{n}W_{n,\ZZ}
\end{align*}
with $R_n:W_{n+1,\ZZ}\twoheadrightarrow W_{n,\ZZ}$. We denote by $\widehat{W}_\ZZ$ the formal completion of $W_\ZZ$ in the sense of Cartier (see \cite[Section~3.2]{MRT2}). Let $A$ be a ring. Then $\widehat{W}_\ZZ$ is characterized as
\begin{align}\label{eq:orig-2-1}
\widehat{W}_\ZZ(A)=\left\{(x_0,x_1,x_2,\ldots)\in W_{\ZZ}(A)\ \middle|\ 
\begin{lgathered}\textrm{$x_i$ is nilpotent for all $i$ and}\\\textrm{$x_i=0$ for all but finitely many $i$}
\end{lgathered}\right\}.
\end{align}

\subsection{Some morphisms of Witt vectors}

Let $A$ be a ring and let $B$ be an $A$-algebra. The following identities on phantom components are taken over the universal polynomial ring and then specialized to $A$.
We define a morphism $F:W_A\rightarrow W_A$ by
\begin{align*}
\Phi_i(F(\xx))=\Phi_{i+1}(\xx)\qquad (i\ge 0)
\end{align*}
for $\xx\in W_A(B)$. If $A$ is an $\FF_p$-algebra, $F$ is nothing but the usual Frobenius endomorphism. Let $[\lambda]$ be the Teichm\"uller lift $[\lambda]=(\lambda,0,0,\ldots)\in W(A)$ for $\lambda\in A$. Then we define the endomorphism $F^{(\lambda)}:=F-[\lambda^{p-1}]$ of $W_A$. We define the Verschiebung endomorphism $V:W_A\to W_A$ by
\[
V(\xx):=(0,x_0,x_1,\ldots)
\]
for $\xx=(x_0,x_1,\ldots)\in W_A(B)$.
Sekiguchi and Suwa \cite[Section~4, p.~20]{SS2} introduced an important morphism called the $T$-map as follows. For $\aa=(a_0,a_1,\ldots)\in W(A)$, we define a morphism $T_{\aa}:W_A\rightarrow W_A$ by
\begin{align*}
\Phi_n(T_{\aa}(\xx))
={a_0}^{p^n}\Phi_n(\xx)+p{a_1}^{p^{n-1}}\Phi_{n-1}(\xx)+\cdots+p^na_n\Phi_0(\xx)
\qquad (n\ge 0)
\end{align*}
for $\xx\in W_A(B)$. 

\section{Deformed Artin-Hasse exponential series}

In this section, we recall several formal series and relations between them, which are introduced and proved in \cite{SS1} and \cite{SS2}.

\subsection{The group scheme $\gg^{(\lambda)}$}

Let $A$ be a ring and let $\lambda\in A$. Put
\begin{align*}
\gg^{(\lambda)}_A:=\Spec\,A\left[X,\dfrac{1}{1+\lambda X} \right].
\end{align*}
We define a morphism
\begin{align*}
\alpha^{(\lambda)}:\gg^{(\lambda)}_A\rightarrow\GG_{m,A};\ x\mapsto 1+\lambda x.
\end{align*}
We endow $\gg^{(\lambda)}_A$ with the commutative group scheme structure whose group law is
\begin{align*}
x\cdot y:=x+y+\lambda xy.
\end{align*}
Then $\alpha^{(\lambda)}$ is a homomorphism of group schemes.
If $\lambda$ is invertible in $A$, then $\alpha^{(\lambda)}$ is an isomorphism over $A$. On the other hand, if $\lambda=0$, then $\gg^{(\lambda)}$ is the additive group scheme $\GG_{a,A}$.

\subsection{Deformed Artin-Hasse exponential series}

The Artin-Hasse exponential series $E_p(X)$ is given by
\begin{align*}
E_p(X)=\exp\left(\sum_{r \geq 0}\frac{X^{p^r}}{p^r}\right).
\end{align*}
The Artin-Hasse exponential series $E_p(X)$ lies in $\ZZ_{(p)}[[X]]$ (\cite[Chap.~V, \S16]{SAG}). We define a formal power series $E_p(U,\Lambda;X)$ in $\QQ[U,\Lambda][[X]]$ by
\begin{align*}
E_p(U,\Lambda;X)=
(1+\Lambda X)^{\frac{U}{\Lambda}}\prod_{k=1}^{\infty}(1+\Lambda^{p^k}X^{p^k})^{\frac{1}{p^k}\left\{\left(\frac{U}{\Lambda}\right)^{p^k}-\left(\frac{U}{\Lambda}\right)^{p^{k-1}}\right\}}.
\end{align*}
As in \cite[Corollary~2.5]{SS1} or \cite[Lemma~4.8]{SS2}, we see that the formal power series $E_p(U,\Lambda;X)$ lies in $\ZZ_{(p)}[U,\Lambda][[X]]$. Note that $E_p(1,0;X)=E_p(X)$. Let $A$ be a $\ZZ_{(p)}$-algebra. For $\lambda\in A$ and $\vv=(v_0,v_1,\ldots)\in W(A)$, we define a formal power series $E_p(\vv,\lambda;X)$ in $A[[X]]$ by
\begin{align*}
E_p(\vv,\lambda;X)&=\prod_{k=0}^{\infty}E_p(v_k,\lambda^{{p^k}};X^{p^k})\\
                  &=(1+\lambda X)^{\frac{v_0}{\lambda}}\prod_{k=1}^{\infty}(1+\lambda^{p^k}X^{p^k})^{\frac{1}{p^k\lambda^{p^k}}\Phi_{k-1}(F^{(\lambda)}(\vv))}.
\end{align*}
If $\lambda$ is nilpotent and $\vv\in\widehat{W}(A)$, by \cite[Proposition~2.11]{SS1} we have $E_p(\vv,\lambda;X)\in A[X]$. Moreover, we define a formal power series $F_p(\vv,\lambda;X,Y)$ as follows:
\begin{align*}
F_p(\vv,\lambda;X,Y)=\prod_{k=1}^{\infty}\left(\frac{(1+\lambda^{p^k}X^{p^k})(1+\lambda^{p^k}Y^{p^k})}{1+\lambda^{p^k}(X+Y+\lambda XY)^{p^k}}\right)^{\frac{1}{p^k\lambda^{p^k}}\Phi_{k-1}(\vv)}.
\end{align*}
The coefficients of $F_p(\vv,\lambda;X,Y)$ belong to $A$; see \cite[Lemma~2.16]{SS1} or \cite[Lemma~4.9]{SS2}.

\subsection{The group scheme $\EE_A^{(\lambda_1,\lambda_2;D)}$ which is an extension of $\gg^{(\lambda_1)}_A$ by $\gg^{(\lambda_2)}_A$}\label{subsec:two-dimensional-exponentials}

Let $A$ be a $\ZZ_{(p)}$-algebra and let $\lambda_1,\lambda_2\in A$ be nonzero divisors. Put $A_{\lambda_2}:=A/(\lambda_2)$. Assume $\lambda_1$ is nilpotent in $A_{\lambda_2}$.
We set
\begin{align*}
\EE^{(\lambda_1,\lambda_2;D)}_A:=\Spec\,A\left[X,Y,\dfrac{1}{1+\lambda_1 X},\dfrac{1}{D(X)+\lambda_2 Y}\right],
\end{align*}
where $\uu\in\widehat{W}(A_{\lambda_2})^{F^{(\lambda_1)}}$ and $D(X)\in A[X]$ is a polynomial lift of $E_p(\uu,\lambda_1;X)$ normalized by $D(0)=1$. We define a morphism $\alpha^{(\lambda_1,\lambda_2)}$ by
\begin{align*}
\alpha^{(\lambda_1,\lambda_2)}:\EE^{(\lambda_1,\lambda_2;D)}_A\rightarrow \GG_{m,A}^2;\ 
(x,y)\mapsto (1+\lambda_1 x,\; D(x)+\lambda_2 y).
\end{align*}
It is known that $\EE^{(\lambda_1,\lambda_2;D)}_A$ admits a unique group scheme structure over $A$ such that $\alpha^{(\lambda_1,\lambda_2)}$ is a homomorphism of group schemes. If $\lambda_1$ and $\lambda_2$ are units in $A$, then $\alpha^{(\lambda_1,\lambda_2)}$ is an isomorphism, and hence $\EE^{(\lambda_1,\lambda_2;D)}\simeq \GG_{m,A}^2$. After a base change on which both parameters vanish, the group is an extension of $\GG_a$ by $\GG_a$; this extension need not split.

For the following formal identities, choose a lift $\uu^\#\in W(A)$ of $\uu$ and put
\begin{align*}
D=D(X):=E_p(\uu^\#,\lambda_1;X).
\end{align*}
If $D_0$ is the polynomial lift above, the coordinate change $(X,Y)\mapsto(X,Y+(D_0-D)/\lambda_2)$ identifies the two formal completions; see also the beginning of Section~5. The matrix $U_2$ uses $\bb=F^{(\lambda_1)}(\uu^\#)/\lambda_2$. We define a formal power series $E_p(\vv,(\lambda_1,\lambda_2);X,Y)$ by
\begin{align*}
E_p(\vv,(\lambda_1,\lambda_2);X,Y)
=E_p(\vv_1,\lambda_1;X)\cdot E_p\!\left(\vv_2,\lambda_2;\frac{Y}{D(X)}\right)\in A[[X,Y]]
\end{align*}
for $\vv=(\vv_1,\vv_2)\in W^2(A)$.

\begin{rem}\label{rem:exp-inj}
{\rm
We will use injectivity of the exponential maps on all Witt vectors, before restricting to a Frobenius kernel. For every $A$-algebra $C$ and any $\lambda\in C$, the map
\begin{align*}
W(C)\rightarrow 1+TC[[T]];\ \zz\mapsto E_p(\zz,\lambda;T)
\end{align*}
is injective. Indeed, recall that
\begin{align*}
E_p(\zz,\lambda;T)=\prod_{k\geq0}E_p(z_k,\lambda^{p^k};T^{p^k}),
\end{align*}
and that $E_p(U,\Lambda;X)=1+UX+O(X^2)$. Hence we have
\begin{align*}
E_p(z_k,\lambda^{p^k};T^{p^k})
=1+z_kT^{p^k}+O(T^{2p^k}).
\end{align*}
For a fixed $n\geq0$, the coefficient of $T^{p^n}$ in $E_p(\zz,\lambda;T)$ receives no contribution from the factors with $k>n$, while the factor with $k=n$ contributes exactly $z_n$. All the remaining contributions depend only on
$z_0,\ldots,z_{n-1}$ and $\lambda$. Therefore, we obtain
\begin{align*}
[T^{p^n}]E_p(\zz,\lambda;T)
=z_n+P_n(z_0,\ldots,z_{n-1},\lambda),
\end{align*}
where $P_n\in\ZZ_{(p)}[Z_0,\ldots,Z_{n-1},\Lambda]$ and $P_0=0$. Here, as usual, $[T^m]f(T)=a_m$ for $f(T)=1+a_1T+a_2T^2+\cdots$. It follows that the Witt coordinates can be recovered recursively from
the coefficients of $E_p(\zz,\lambda;T)$. Namely,
\begin{align*}
z_0=[T]E_p(\zz,\lambda;T)\quad \text{and}\quad
z_n=[T^{p^n}]E_p(\zz,\lambda;T)-P_n(z_0,\ldots,z_{n-1},\lambda) \quad(n\geq1).
\end{align*}
Thus, if $E_p(\zz,\lambda;T)=E_p(\zz',\lambda;T)$, then $z_0=z_0'$, and induction on $n$ gives $z_n=z_n'$ for every $n\geq0$. Hence $\zz=\zz'$, proving the injectivity. The two-dimensional map
\begin{align*}
W^2(C)\rightarrow 1+(X,Y)C[[X,Y]];\ \vv\longmapsto E_p(\vv,(\lambda_1,\lambda_2);X,Y)
\end{align*}
is also injective: setting $Y=0$ recovers the first factor, and setting $X=0$ recovers the second factor, since $D(0)=1$. These arguments apply after every base change, including when the parameters are zero divisors.
}
\end{rem}

Let $\XX=(X_0,X_1)$ and $\YY=(Y_0,Y_1)$. We also consider a formal power series
\begin{align*}
F_p(\vv,(\lambda_1,\lambda_2);\XX,\YY)\in A[[\XX,\YY]]
\end{align*}
for $\vv\in W^2(A)$ (see \cite[p.~34]{SS2} for details). Note that, for $\vv\in W^2(A)$, the equality
\begin{align}\label{eq:orig-3-1}
F_p(U_2(\vv),(\lambda_1,\lambda_2);\XX,\YY)
=\frac{E_p(\vv,(\lambda_1,\lambda_2);\XX)\cdot E_p(\vv,(\lambda_1,\lambda_2);\YY)}{E_p(\vv,(\lambda_1,\lambda_2);\XX\cdot\YY)}
\end{align}
holds, where $\XX\cdot\YY$ denotes the product in the formal group $\widehat{\EE}^{(\lambda_1,\lambda_2;D)}_A$ (see \cite[Theorem~4.16]{SS2}). The formal-group isomorphisms of \cite[Theorem~5.1]{SS2} take the following form:
\begin{align}\label{eq:orig-3-2}
\begin{split}
  \Ker[U_2:W^2_A(C)\to W^2_A(C)]&\xrightarrow{\sim}\Hom(\widehat{\EE}^{(\lambda_1,\lambda_2;D)}_C,\widehat{\GG}_{m,C});\\
                             \vv&\mapsto E_p(\vv,(\lambda_1,\lambda_2);X,Y)
\end{split}
\end{align}
and
\begin{align}\label{eq:orig-3-3}
\begin{split}
\Coker[U_2:W^2_A(C)\to W^2_A(C)]&\xrightarrow{\sim} H^2_0(\widehat{\EE}^{(\lambda_1,\lambda_2;D)}_C,\widehat{\GG}_{m,C});\\
                             \bigl[\vv\bigr]&\mapsto [F_p(\vv,(\lambda_1,\lambda_2);\XX,\YY)]
\end{split}
\end{align}
for any $A$-algebra $C$.

\begin{rem}
{\rm
We give the comparison in \eqref{eq:orig-3-3} over an arbitrary
$A$-algebra $C$. For \eqref{eq:orig-3-2}, see also
\cite[Theorem~4.1(2)]{MRT2}. Put
\[
G=\widehat{\EE}^{(\lambda_1,\lambda_2;D)}_C,\qquad
G_i=\widehat{\gg}^{(\lambda_i)}_C,\qquad H=\widehat{\GG}_{m,C}.
\]
The projection $\pi:G\to G_1$ has the section $s(x)=(x,0)$ as a
map of pointed formal schemes. Write $i:G_2\to G$, $y\mapsto(0,y)$.
Since $D(0)=1$, we have $(x,t)=i(t/D(x))s(x)$. Define
$h(x,x')$ by $s(x)s(x')=i(h(x,x'))s(xx')$, and put
$d_H(\chi)=[\chi(h)^{-1}]$ for $\chi\in\Hom(G_2,H)$.

We first check exactness at the two middle terms of
\[
\Hom(G_2,H)\xrightarrow{d_H} H^2_0(G_1,H)
\xrightarrow{\pi^*}H^2_0(G,H)\xrightarrow{i^*}H^2_0(G_2,H).
\]
Here $(\delta a)(g,g')=a(g)a(g')/a(gg')$.
If $\pi^*c=\delta a$, then $\chi=a\circ i$ is a character and
$a(i(y)s(x))=\chi(y)a(s(x))$. Hence
$c=\delta(a\circ s)\,\chi(h)^{-1}$. Conversely,
$a(i(y)s(x))=\chi(y)$ gives $\delta a=\pi^*(\chi(h)^{-1})$.
Next, let $c$ be a normalized symmetric cocycle on $G$ with
$i^*c=\delta b$. Put
\[
a(i(y)s(x))=\frac{b(y)}{c(i(y),s(x))},\qquad c'=c/\delta a.
\]
The cocycle identity on $(i(k),i(y),s(x))$ gives
$c'(i(k),i(y)s(x))=1$. By the cocycle identity and symmetry,
$c'(i(k)g,g')=c'(g,g')=c'(g,i(k)g')$.
Thus $c'=\pi^*((s\times s)^*c')$, proving exactness.

Put $f_i=F^{(\lambda_i)}$, $K_i=\Ker(f_i(C))$,
$P_i=\Coker(f_i(C))$ and $d:K_2\to P_1$, $\tau\mapsto[-T_{\bb}\tau]$.
By the definition of $F_p$ in \cite[(45), p.~34]{SS2}, we have
\[
F_p((\xi,0))=\pi^*F_p(\xi,\lambda_1),\qquad
i^*F_p((\xi_1,\xi_2))=F_p(\xi_2,\lambda_2),
\]
where the variables are omitted. The second equality and
\cite[Theorem~2.19.1]{SS1} give the surjectivity of $i^*$ on $H^2_0$.
For $\tau\in K_2$, the cochain $a=E_p((0,\tau),(\lambda_1,\lambda_2))$
satisfies, by \eqref{eq:orig-3-1},
\[
a\circ i=E_p(\tau,\lambda_2),\qquad
\delta a=\pi^*F_p(-T_{\bb}\tau,\lambda_1).
\]
Hence $d_H(E_p(\tau,\lambda_2))=[F_p(-T_{\bb}\tau,\lambda_1)]$.
By \eqref{eq:orig-3-1}, the map $\vv\mapsto[F_p(\vv)]$ factors through $\Coker(U_2(C))$.
The above equalities give a map between the exact sequences
\[
\begin{aligned}
0&\longrightarrow P_1/\Im d\longrightarrow\Coker(U_2(C))
  \longrightarrow P_2\longrightarrow0,\\
0&\longrightarrow H^2_0(G_1,H)/\Im d_H\longrightarrow H^2_0(G,H)
  \longrightarrow H^2_0(G_2,H)\longrightarrow0,
\end{aligned}
\]
where the first sequence follows from the triangular form of $U_2$.
The two outer maps are isomorphisms by \cite[Theorem~2.19.1]{SS1}.
Therefore the middle map is an isomorphism, as asserted.
}
\end{rem}

\section{The $l$-th Kummer--Artin--Schreier--Witt type isogeny}

In this section, we construct the two-dimensional Kummer--Artin--Schreier--Witt type isogeny.

\subsection{Review: a discrete valuation ring containing a root of unity}

Let $\Phi_n(X)$ be the $n$-th cyclotomic polynomial. Let $m$ be a positive integer. By the equality $X^n-1=\prod_{d\mid n}\Phi_d(X)$, we have
\begin{align*}
\Phi_{p^m}(X)
&=\dfrac{\prod_{d\mid p^m}\Phi_d(X)}{\prod_{d\mid p^{m-1}}\Phi_d(X)}
=\dfrac{X^{p^m}-1}{X^{p^{m-1}}-1}\\
&=X^{p^{m-1}(p-1)}+X^{p^{m-1}(p-2)}+\cdots+X^{p^{m-1}}+1.
\end{align*}
In particular, we obtain $\Phi_{p^m}(1)=p$. Put $f(X):=\Phi_{p^m}(X+1)$. Then
\begin{align*}
f(X)\equiv X^{p^m-p^{m-1}}\ (\mod\ p)\quad \textrm{and}\quad f(0)=\Phi_{p^m}(1)=p
\end{align*}
hold. Hence $f(X)=\Phi_{p^m}(X+1)$ is an Eisenstein polynomial. Therefore $\ZZ_{(p)}[X]/(f(X))$ is a discrete valuation ring with uniformizer $\overline{X}$ (\cite[Chapter~I, Proposition~17]{SL}). Let $\zeta_{p^m}$ be a primitive $p^m$-th root of unity and set $\lambda:=\zeta_{p^m}-1$. For $f(X)=\Phi_{p^m}(X+1)$, we have $f(\lambda)=\Phi_{p^m}(\zeta_{p^m})=0$. Hence the isomorphism
\begin{align*}
\ZZ_{(p)}[X]/(f(X))\xrightarrow{\sim}\ZZ_{(p)}[\lambda]=\ZZ_{(p)}[\zeta_{p^m}];\ \overline{X}\mapsto\lambda
\end{align*}
holds. Note that, for the discrete valuation ring $A:=\ZZ_{(p)}[\zeta_{p^m}]$ with uniformizer $\lambda$, we obtain $\mathrm{Frac}(A)=\QQ(\zeta_{p^m})$ and $A/(\lambda)=\FF_p$. Let $v$ be the valuation of $A$, normalized by $v(\lambda)=1$. Since $v(p)=\varphi(p^m)=p^m-p^{m-1}$, we have
\begin{align}\label{eq:orig-4-1}
p=u\lambda^{p^m-p^{m-1}}\quad (u\in A^\times).
\end{align}

\begin{rem}\label{st:orig-4-1}
{\rm
Assume $l\leq m$. Take positive integers $m_1$ and $m_2$ such that $l\leq m_i\leq m$ for $i=1,2$. Put $\lambda_i:=\zeta_{p^{m_i}}-1\in A$ for $i=1,2$. By \eqref{eq:orig-4-1}, we have
\begin{align*}
p=u_i\lambda_i^{p^{m_i}-p^{m_i-1}}\quad (u_i\in A^\times)
\end{align*}
for $i=1,2$. For each $0\leq k<l$, we put
\begin{align*}
\nu_{i,k}:=u_i^{l-k}\lambda_i^{(l-k)(p^{m_i}-p^{m_i-1})-(p^l-p^k)}.
\end{align*}
Since
\begin{align*}
p^l-p^k=(p-1)\sum_{i=0}^{l-k-1}p^{k+i}\leq (p-1)(l-k)p^{l-1}=(l-k)(p^l-p^{l-1}),
\end{align*}
we have $\nu_{i,k}\in A$. Note that, by \eqref{eq:orig-4-1}, if $l=m_1=m_2$, then, for each $i=1,2$, $\nu_{i,k}$ is divisible by $\lambda_i$ for $0\leq k<l-1$, and $\nu_{i,l-1}$ is a unit. If $l<m_i\leq m$, then $\nu_{i,k}$ is divisible by $\lambda_i$ for $0\leq k<l$.
}
\end{rem}

\subsection{The one-dimensional formal case}

We establish the identities needed in Section~5. The corresponding pullbacks were studied in \cite{A4}. We give proofs that account for the correction factor for arbitrary Witt vectors and for the hypotheses needed for formal surjectivity. Temporarily, omit the subscript $i$ and write $\lambda$, $\nu_k$, and $\aa$. Here $\lambda$ is a nonzero divisor in $A$; its images after base change need not be.

\begin{rem}\label{st:orig-4-9}
{\rm
Write
\begin{align*}
p^l[1]=(c_0,c_1,c_2,\ldots)\in W(\ZZ).
\end{align*}
An induction using the phantom components shows that $c_k\in p^{l-k}\ZZ$ for $0\leq k\leq l-1$.  Write $c_k=p^{l-k}d_k$ with $d_k\in\ZZ$ for these values of $k$. Since $p^l[\Lambda]=[\Lambda]\,p^l[1]$, multiplication by the Teichm\"uller vector gives
\begin{align*}
p^l[\Lambda]
=\left(c_0\Lambda,c_1\Lambda^p,c_2\Lambda^{p^2},\ldots\right).
\end{align*}
For $0\leq k\leq l-1$, the $k$-th component of $p^l[\Lambda]/\Lambda^{p^l}$ is $\left(d_kp^{l-k}\right)/\Lambda^{p^l-p^k}$, whereas for $k\geq l$ it is $c_k\Lambda^{p^k-p^l}$.  Consequently,
\begin{align*}
\frac{1}{\Lambda^{p^l}}p^l[\Lambda]
\in
W\left(\ZZ\left[\Lambda,\left\{\frac{p^{l-k}}{\Lambda^{p^l-p^k}}\ \middle|\ 0\leq k\leq l-1\right\}\right]\right).
\end{align*}
By \eqref{eq:orig-1-1}, we may specialize $\Lambda$ to $\lambda$ and $p^{l-k}/\Lambda^{p^l-p^k}$ to $\nu_k$ for $0\leq k\leq l-1$. This shows that every component of $p^l[\lambda]/\lambda^{p^l}$ lies in $A$.
}
\end{rem}

Assume that
\begin{align}\label{eq:orig-4-2}
p^{l-k}=\nu_k\lambda^{p^l-p^k}\quad(0\leq k\leq l-1).
\end{align}
We put
\begin{align*}
\psi^{(l)}(X)
:=\dfrac{1}{\lambda^{p^l}}\left\{(1+\lambda X)^{p^l}-1\right\}
=X^{p^l}+\sum_{r=1}^{p^l-1}\alpha_rX^r.
\end{align*}
For $1\leq r\leq p^l-1$, set $k_r:=v_p(r)$. Then $\binom{p^l}{r}=p^{l-k_r}c_r$ holds for some $c_r\in\ZZ$. Hence, by \eqref{eq:orig-4-2}, we have
\begin{align}\label{eq:orig-4-3}
\alpha_r=\binom{p^l}{r}\lambda^{r-p^l}=c_r\nu_{k_r}\lambda^{r-p^{k_r}}.
\end{align}
In particular, if the image of every $\nu_k$ is nilpotent in an $A$-algebra $B$, then the image of every $\alpha_r$ is nilpotent in $B$. Assume that $B$ is an $A$-algebra such that each $\nu_k$ is nilpotent.

\begin{lem}\label{st:orig-4-2}
The homomorphism
\begin{align*}
\widehat{\psi}^{(l)}:\widehat{\mathcal G}^{(\lambda)}_B\rightarrow\widehat{\mathcal G}^{(\lambda^{p^l})}_B
\end{align*}
induced by $\psi^{(l)}$ is an epimorphism for the fppf topology.
\end{lem}

\begin{proof}
Let $C$ be a $B$-algebra and let $y\in\widehat{\mathcal G}^{(\lambda^{p^l})}_B(C)$. Since $\widehat{\mathcal G}^{(\lambda^{p^l})}_B$ is the formal completion along the zero section, $y$ is nilpotent in $C$. Consider
\begin{align*}
C':=
C[T]/\left(\psi^{(l)}(T)-y\right).
\end{align*}
Since $\psi^{(l)}(T)-y$ is monic of degree $p^l$, the algebra $C'$ is finite free of rank $p^l$ over $C$. Hence $\Spec\,C'\rightarrow\Spec\,C$ is an fppf covering. Let $t$ be the image of $T$ in $C'$. By \eqref{eq:orig-4-3}, all the elements $\alpha_1,\ldots,\alpha_{p^l-1}$ are nilpotent in $C'$. Put
\begin{align*}
J:=(y,\alpha_1,\ldots,\alpha_{p^l-1})C'.
\end{align*}
Since $J$ is generated by finitely many nilpotent elements, it is a
nilpotent ideal. The relation $\psi^{(l)}(t)=y$ gives
\begin{align*}
t^{p^l}
=y-\sum_{r=1}^{p^l-1}\alpha_rt^r\in J.
\end{align*}
Since $J$ is nilpotent, $t$ is nilpotent. Therefore $t\in\widehat{\mathcal G}^{(\lambda)}_B(C')$ and $\widehat{\psi}^{(l)}_B(t)=y$. Thus $\widehat{\psi}^{(l)}$ is an epimorphism for the fppf topology.
\end{proof}

As a consequence, we obtain the injectivity of pullback on characters.

\begin{lem}\label{st:orig-4-3}
For every $B$-algebra $C$, the pullback homomorphism
\begin{align*}
\left(\widehat{\psi}^{(l)}\right)^\ast:
\Hom\!\left(\widehat{\mathcal G}_C^{(\lambda^{p^l})},\widehat{\GG}_{m,C}\right)
\rightarrow
\Hom\!\left(\widehat{\mathcal G}^{(\lambda)}_C,\widehat{\GG}_{m,C}\right)
\end{align*}
is injective.
\end{lem}

\begin{proof}
This follows immediately from Lemma~\ref{st:orig-4-2}.
\end{proof}

For a Witt vector $\zz=(z_0,z_1,\ldots)$, put
\begin{align*}
\widetilde p^{\,r}\zz:=V^r(z_0^{p^r},z_1^{p^r},\ldots)\quad(r\geq0).
\end{align*}
We define
\begin{align*}
\widetilde p^{\,r}E_p(\zz,\lambda;T):=E_p(\widetilde p^{\,r}\zz,\lambda;T).
\end{align*}
Note that even if
\begin{align*}
E_p(\zz,\lambda;T)=E_p(\zz',\lambda';f(T)),
\end{align*}
it does not necessarily follow that
\begin{align*}
E_p(\widetilde p^{\,r}\zz,\lambda;T)=E_p(\widetilde p^{\,r}\zz',\lambda';f(T)).
\end{align*}
Put
\begin{align*}
E:=E_p(p^l[\lambda],\lambda;X)=(1+\lambda X)^{p^l}.
\end{align*}
We retain the notation $G_p$ of \cite[(34), p.~29]{SS2}. For a Witt vector of variables $\boldsymbol\eta$, set
\begin{align}\label{eq:one-dimensional-correction}
G_p(\boldsymbol\eta,\lambda^{p^l};E)
:=\prod_{r\geq1}\left(\frac{1+(E-1)^{p^r}}{\widetilde p^{\,r}E}\right)^{\frac{1}{p^r\lambda^{p^{l+r}}}\Phi_{r-1}(\boldsymbol\eta)}.
\end{align}
Here $E$ retains the displayed exponential presentation, so $\widetilde p^{\,r}E=E_p(\widetilde p^{\,r}(p^l[\lambda]),\lambda;X)$.

\begin{lem}\label{st:orig-4-4}
For every $A$-algebra $C$ and every $\vv\in W(C)$, we have
\begin{align}\label{eq:one-dimensional-pullback}
\begin{split}
E_p(\vv,\lambda^{p^l};\psi^{(l)}(X))
=E_p(T_{\aa}(\vv),\lambda;X)G_p(F^{(\lambda^{p^l})}(\vv),\lambda^{p^l};E).
\end{split}
\end{align}
In particular, if $F^{(\lambda^{p^l})}(\vv)=0$, then
$E_p(\vv,\lambda^{p^l};\psi^{(l)}(X))=E_p(T_{\aa}\vv,\lambda;X)$.
\end{lem}

\begin{proof}
By \cite[Proposition~4.12]{SS2}, the coefficients of
$G_p(\boldsymbol\eta,\Lambda_2;E_p(\mathbf U,\Lambda_1;X))$
belong to
\begin{align*}
\ZZ_{(p)}[\boldsymbol\eta,\mathbf{U}/\Lambda_2,\Lambda_1,\Lambda_2],
\end{align*}
where $\mathbf U/\Lambda_2$ denotes componentwise division. Specialize by
\begin{align*}
\Lambda_1\mapsto\lambda,\quad
\Lambda_2\mapsto\lambda^{p^l}\quad\text{and}\quad
\mathbf U/\Lambda_2\mapsto\aa.
\end{align*}
Since $\aa=p^l[\lambda]/\lambda^{p^l}\in W(A)$, this is a well-defined specialization and sends $\mathbf U$ to $p^l[\lambda]$. Thus \eqref{eq:one-dimensional-correction} belongs to $1+XA[\boldsymbol\eta][[X]]$.

We follow the calculation in \cite[proof of Lemma~1, p.~123]{A1}, retaining the correction factor since $\vv$ is not assumed to lie in $\Ker(F^{(\lambda^{p^l})})(C)$. By \cite[Proposition~4.11]{SS2}
and \cite[(36), p.~29]{SS2}, we obtain
\begin{align*}
E_p(\vv,\lambda^{p^l};\psi^{(l)}(X))
&=E_p\!\left(\vv,\lambda^{p^l};\frac{E-1}{\lambda^{p^l}}\right)\\
&=\widetilde E_p(\vv,\lambda^{p^l};E)\,
G_p(F^{(\lambda^{p^l})}(\vv),\lambda^{p^l};E)\\
&=E_p(T_{\aa}(\vv),\lambda;X)\,
G_p(F^{(\lambda^{p^l})}(\vv),\lambda^{p^l};E).
\end{align*}
Finally, $G_p(0,\lambda^{p^l};E)=1$ by its definition.
\end{proof}

\begin{lem}\label{st:orig-4-5}
There is a unique group scheme homomorphism $T'_{\aa}:W_A\to W_A$ such that
\begin{align}\label{eq:one-dimensional-factorization}
F^{(\lambda)}\circ T_{\aa}=T'_{\aa}\circ F^{(\lambda^{p^l})}.
\end{align}
If every $\nu_k$ is nilpotent in an $A$-algebra $B$, then $T_{\aa,B}$ and $T'_{\aa,B}$ are closed immersions. In particular, for every $B$-algebra $C$,
\begin{align}\label{eq:orig-4-4}
\Ker(F^{(\lambda^{p^l})})(C)=\Ker(F^{(\lambda)}\circ T_{\aa})(C).
\end{align}
\end{lem}

\begin{proof}
Write
\begin{align*}
f:=F^{(\lambda^{p^l})}:W_A\rightarrow W_A\quad\text{and}\quad
h:=F^{(\lambda)}\circ T_{\aa}:W_A\rightarrow W_A.
\end{align*}
We first construct $T'_{\aa}$ over $A$. We then use the nilpotency
assumption over $B$ to prove that $T_{\aa,B}$ and $T'_{\aa,B}$ are
closed immersions. The equality of kernels will follow from the
injectivity of $T'_{\aa,B}$.

\smallskip
\noindent\emph{Step 1: $h$ vanishes on $\Ker(f)$.}
The character description of \cite[Theorem~2.19.1]{SS1} gives, for every $A$-algebra $C$,
\begin{align}\label{eq:orig-4-5}
\Ker(F^{(\lambda)})(C)
\xrightarrow{\sim}
\Hom(\widehat{\mathcal G}^{(\lambda)}_C,\widehat{\GG}_{m,C});\ 
\zz\mapsto E_p(\zz,\lambda;X).
\end{align}
The same statement holds with $\lambda$ replaced by $\lambda^{p^l}$. Take $\vv\in\Ker(f)(C)$. Thus $F^{(\lambda^{p^l})}(\vv)=0$, and
$E_p(\vv,\lambda^{p^l};X)$ is a character of $\widehat{\mathcal G}^{(\lambda^{p^l})}_C$. Pulling this character back along $\widehat\psi^{(l)}_C$ gives a character of $\widehat{\mathcal G}^{(\lambda)}_C$. By Lemma~\ref{st:orig-4-4},
its defining series is
\begin{align*}
E_p(\vv,\lambda^{p^l};\psi^{(l)}(X))=E_p(T_{\aa}(\vv),\lambda;X).
\end{align*}
By the surjectivity of \eqref{eq:orig-4-5}, there is therefore $\zz\in\Ker(F^{(\lambda)})(C)$ such that
\begin{align*}
E_p(\zz,\lambda;X)=E_p(T_{\aa}(\vv),\lambda;X).
\end{align*}
The exponential is injective on all of $W(C)$ by Remark~\ref{rem:exp-inj}, so $\zz=T_{\aa}(\vv)$. Consequently,
\begin{align*}
h(\vv)=F^{(\lambda)}\circ T_{\aa}(\vv)=F^{(\lambda)}(\zz)=0.
\end{align*}
Since $C$ and $\vv$ were arbitrary, this proves that the restriction of $h$ to $\Ker(f)$ is the zero morphism of group schemes.

\smallskip
\noindent\emph{Step 2: construction and uniqueness of $T'_{\aa}$.}
The morphism $f=F^{(\lambda^{p^l})}$ is faithfully flat by
\cite[Corollary~1.8]{SS1}; see also \cite[Lemma~3.7]{MRT2}.
It is affine and hence quasi-compact, so it is an epimorphism for
the fpqc topology. We explain explicitly how Step~1 allows $h$ to
descend along $f$.

Let $C$ be an $A$-algebra and take $\ww\in W(C)$. After a faithfully
flat base change $C\to C'$, choose $\vv\in W(C')$ with
\begin{align*}
f(\vv)=\ww|_{C'}.
\end{align*}
Such a choice exists by taking the base change of $f$ along the
point $\ww:\Spec C\to W_A$.
Let $\vv_1,\vv_2\in W(C'\otimes_C C')$ be the two images of
$\vv$. They have the same image under $f$, so
\begin{align*}
\vv_1-\vv_2\in\Ker(f)(C'\otimes_C C').
\end{align*}
Step~1 and additivity of $h$ imply
\begin{align*}
h(\vv_1)-h(\vv_2)=h(\vv_1-\vv_2)=0.
\end{align*}
Since the representable functor $W_A$ is an fpqc sheaf,
$h(\vv)\in W(C')$ descends to a unique element of $W(C)$.
Define this element to be $T'_{\aa}(\ww)$.

This definition is independent of the chosen cover and lift:
on a common faithfully flat cover, any two lifts differ by an
element of $\Ker(f)$, on which $h$ vanishes. The same construction
shows compatibility with base change. It is additive, since sums
of lifts lift sums and $h$ is additive. Thus we obtain a natural
homomorphism $W(C)\to W(C)$ for every $A$-algebra $C$.
By Yoneda, these maps define a group scheme homomorphism
$T'_{\aa}:W_A\to W_A$, and the construction gives
\begin{align*}
T'_{\aa}\circ f=h.
\end{align*}
This is \eqref{eq:one-dimensional-factorization}.
If another morphism has the same composite with $f$, its value at
any $\ww$ agrees with $T'_{\aa}(\ww)$ after the above cover, and
hence already over $C$. This proves uniqueness.

\smallskip
\noindent\emph{Step 3: reduction of both morphisms to $V^l$.}
Now let $B$ satisfy the nilpotency assumption, and put
\begin{align*}
I=(\nu_0,\ldots,\nu_{l-1})B,
\qquad \overline B=B/I.
\end{align*}
The ideal $I$ is nilpotent because it is generated by finitely many
nilpotent elements. Indeed, if $\nu_k^{e_k}=0$ for each $k$, then
$I^N=0$ for $N=1+\sum_{k=0}^{l-1}(e_k-1)$.
Also,
\begin{align*}
p=\nu_{l-1}\lambda^{p^l-p^{l-1}}\in I,
\end{align*}
so $\overline B$ is an $\FF_p$-algebra.

We apply the characteristic-$p$ calculation in
\cite[proof of Lemma~1, pp.~123--124]{A1} after identifying
the reduction of the particular vector $\aa$ used here.
Write $p^l[1]=(c_0,c_1,\ldots)\in W(\ZZ)$ and, for $0\leq k<l$,
write $c_k=p^{l-k}d_k$ with $d_k\in\ZZ$, as in
Remark~\ref{st:orig-4-9}. Since
$p^l[\lambda]=(c_k\lambda^{p^k})_{k\geq0}$, the componentwise
division defining $\aa$ over $A$ and \eqref{eq:orig-4-2} give
\begin{align*}
a_k=
\begin{cases}
d_k\nu_k,&0\leq k<l,\\
c_k\lambda^{p^k-p^l},&k\geq l.
\end{cases}
\end{align*}
These identities may be reduced modulo $I$ after base change to
$B$, without dividing by the image of $\lambda$ in $\overline B$.
The calculation in the cited proof, with $\lambda=1$, gives
\begin{align*}
c_k\equiv
\begin{cases}
1\pmod p,&k=l,\\
0\pmod p,&k\ne l.
\end{cases}
\end{align*}
Since $\nu_k\in I$ for $k<l$ and $p\in I$, it follows that
\begin{align*}
\overline{\aa}=V^l[1]\in W(\overline B).
\end{align*}
As in the cited proof, the identity $T_{V^l[1]}=V^l$ now gives
\begin{align*}
\overline{T_{\aa}}=V^l:W_{\overline B}\longrightarrow W_{\overline B}.
\end{align*}

It remains to compute the reduction of $T'_{\aa}$.
All operators in the following calculation are over $\overline B$;
write $\overline\lambda$ for the image of $\lambda$. In characteristic $p$, Frobenius acts componentwise,
so $FV^l=V^lF$. Moreover, Teichm\"uller multiplication satisfies
\begin{align*}
[\overline\lambda^{p-1}]\circ V^l
=V^l\circ[\overline\lambda^{p^l(p-1)}].
\end{align*}
Consequently,
\begin{align*}
F^{(\overline\lambda)}\circ V^l
&=(F-[\overline\lambda^{p-1}])\circ V^l\\
&=V^l\circ(F-[\overline\lambda^{p^l(p-1)}])
=V^l\circ F^{(\overline\lambda^{p^l})}.
\end{align*}
Reducing \eqref{eq:one-dimensional-factorization} modulo $I$
therefore yields
\begin{align*}
\overline{T'_{\aa}}\circ F^{(\overline\lambda^{p^l})}
=V^l\circ F^{(\overline\lambda^{p^l})}.
\end{align*}
Faithful flatness is preserved by base change, so
$F^{(\overline\lambda^{p^l})}$ is still an fpqc epimorphism. We may cancel
it on the right, obtaining
\begin{align*}
\overline{T'_{\aa}}=V^l.
\end{align*}

\smallskip
\noindent\emph{Step 4: lifting the closed immersions and comparing kernels.}
Let $\theta$ be either $T_{\aa,B}$ or $T'_{\aa,B}$. Write its
coordinate ring map as
\begin{align*}
\theta^*:B[Z_0,Z_1,\ldots]\longrightarrow B[T_0,T_1,\ldots].
\end{align*}
By Step~3, its reduction modulo $I$ is the coordinate map of $V^l$,
given by
\begin{align*}
Z_n\longmapsto
\begin{cases}
0,&n<l,\\
T_{n-l},&n\geq l.
\end{cases}
\end{align*}
This map is surjective, since every generator $T_m$ is the image
of $Z_{m+l}$.

Let $M=\Coker(\theta^*)$ as a $B$-module. Tensoring the defining
right exact sequence with $B/I$ gives
\begin{align*}
M/IM=\Coker(\theta^*\otimes_B B/I)=0.
\end{align*}
Thus $M=IM$. Since $I^N=0$, iteration gives
\begin{align*}
M=IM=I^2M=\cdots=I^NM=0.
\end{align*}
This argument uses nilpotency of $I$ and does not require $M$ to
be finitely generated. It follows that $\theta^*$ is surjective.
A morphism of affine schemes is a closed immersion exactly when
its coordinate ring map is surjective; hence both $T_{\aa,B}$ and
$T'_{\aa,B}$ are closed immersions.

In particular, $T'_{\aa}(C)$ is injective for every $B$-algebra $C$.
For any $\vv\in W(C)$, \eqref{eq:one-dimensional-factorization}
therefore gives
\begin{align*}
F^{(\lambda)}T_{\aa}\vv=0
&\ \Longleftrightarrow\ T'_{\aa}(F^{(\lambda^{p^l})}\vv)=0\\
&\ \Longleftrightarrow\ F^{(\lambda^{p^l})}\vv=0.
\end{align*}
This proves \eqref{eq:orig-4-4}. Notice that no nonzero-divisor
assumption on the image of $\lambda$ in $B$ or $C$ was used.
\end{proof}

\begin{rem}\label{st:orig-4-6}
{\rm
The general pullback identity requires the factor $G_p(F^{(\lambda^{p^l})}\vv,\lambda^{p^l};E)$.
It becomes $1$ on the Frobenius kernel. In Lemma~\ref{st:orig-4-5}, nilpotency is used to lift the closed-immersion property of $V^l$ from $B/I$ to $B$. Thus the reverse kernel inclusion does not use an uncorrected pullback identity for arbitrary Witt vectors.
}
\end{rem}

\subsection{Proof of Theorem~\ref{st:orig-1-2}}
Let $A$ be a $\ZZ_{(p)}$-algebra, and let $\lambda_1,\lambda_2\in A$ be nonzero divisors satisfying $\lambda_1^{p^l-p^{l-1}}\mid p$ and $\lambda_2^{p^l-p^{l-1}}\mid p$.
\begin{rem}\label{st:orig-4-7}
{\rm
Under the assumption $\lambda_i^{p^l-p^{l-1}}\mid p$ for $i=1,2$, there exists $u_i\in A$ such that $p=u_i\lambda_i^{p^l-p^{l-1}}$. For each $0\leq k<l$, define
\begin{align*}
\nu_{i,k}:=u_i^{l-k}\lambda_i^{(l-k)(p^l-p^{l-1})-(p^l-p^k)}\in A.
\end{align*}
Then we have
\begin{align*}
p^{l-k}=\nu_{i,k}\lambda_i^{p^l-p^k}\quad (0\leq k<l).
\end{align*}
}
\end{rem}

\begin{rem}\label{st:orig-4-8}
{\rm
If $\lambda_1$ and $\lambda_2$ are units in $A$, then \eqref{eq:orig-1-1} holds automatically, and we may take
\begin{align*}
\nu_{1,k}=\frac{p^{l-k}\lambda_1^{p^k}}{\lambda_1^{p^l}}
\quad\text{and}\quad
\nu_{2,k}=\frac{p^{l-k}\lambda_2^{p^k}}{\lambda_2^{p^l}}.
\end{align*}
If $A$ is an $\FF_p$-algebra and we set $\nu_{1,k}=\nu_{2,k}=0$ for all $0\le k\le l-1$, then \eqref{eq:orig-1-1} holds for arbitrary nonzero divisors $\lambda_1,\lambda_2\in A$. If $l\leq m$, a nontrivial example is given in Remark~\ref{st:orig-4-1}.
}
\end{rem}

\begin{sublem}\label{st:orig-4-10}
Let $N$ and $M$ be positive integers such that $N<M$. Then the natural projection
\begin{align*}
\widehat{W}(A_{\lambda_2^{M}})\rightarrow \widehat{W}(A_{\lambda_2^{N}})
\end{align*}
is surjective.
\end{sublem}
\begin{proof}
Take any $\xx=(x_0,x_1,\ldots)\in\widehat{W}(A_{\lambda_2^N})$. By \eqref{eq:orig-2-1}, each $x_i$ is nilpotent in $A_{\lambda_2^N}$, and $x_i=0$ for all but finitely many $i$.  For each nonzero component choose a lift $\widetilde{x_i}\in A$, and choose the zero lift for every zero component. Since $x_i$ is nilpotent, there exists a positive integer $n_i$ such that $x_i^{n_i}=0$ in $A_{\lambda_2^N}$. Hence $\widetilde{x_i}^{\,n_i}=\alpha_i\lambda_2^N$ for some $\alpha_i\in A$. Choose a positive integer $r$ such that $Nr\ge M$. Then
\begin{align*}
\widetilde{x_i}^{\,n_i r}=\alpha_i^r\lambda_2^{Nr}\equiv 0 \pmod{\lambda_2^M}.
\end{align*}
Therefore the image of $\widetilde{x_i}$ in $A_{\lambda_2^M}$ is nilpotent. Since the chosen lifts have finite support, they define an element of $\widehat{W}(A_{\lambda_2^M})$ lifting $\xx$.
\end{proof}
\begin{lem}\label{st:orig-4-11}
Let $D(X)\in A[X]$ be a normalized lift of $E_p(\widetilde{\uu},\lambda_1;X)\in A_{\lambda_2^{p^l}}[X]$, and let $D'(X)\in A[X]$ be a normalized lift of $E_p(\vv,\lambda_1^{p^l};X)\in A_{\lambda_2^{p^l}}[X]$. Then the congruence
\begin{align*}
D(X)^{p^l}(1+\lambda_1 X)^{-j}\equiv D'\left(\psi_1(X)\right)\pmod{\lambda_2^{p^l}}
\end{align*}
holds.
\end{lem}
\begin{proof}
Since $F^{(\lambda_1^{p^l})}(\vv)=0$ over $A/(\lambda_2^{p^l})$, the correction factor in Lemma~\ref{st:orig-4-4} is $1$ over this quotient. Thus \eqref{eq:orig-1-2} and that lemma give
\begin{align*}
D(X)^{p^l}(1+\lambda_1 X)^{-j}&\equiv E_p(\widetilde{\uu},\lambda_1;X)^{p^l}E_p([\lambda_1],\lambda_1;X)^{-j}\\
                              &=E_p(p^l\widetilde{\uu},\lambda_1;X)\cdot E_p(-j[\lambda_1],\lambda_1;X)\\
                              &=E_p(p^l\widetilde{\uu}-j[\lambda_1],\lambda_1;X)\\
                              &\equiv E_p(T_{\aa_1}(\vv),\lambda_1;X)\\
                              &=E_p(\vv,\lambda_1^{p^l};\psi_1(X))\equiv D'(\psi_1(X))\pmod{\lambda_2^{p^l}}
\end{align*}
hold.
\end{proof}
\begin{proof}[Proof of Theorem~\ref{st:orig-1-2}]
The equality
\begin{align*}
\left(\alpha^{(\lambda_1^{p^l},\lambda_2^{p^l})}\circ\Psi_{\uu,\vv}^{(l,j)}\right)(x,y)
=\left(\Theta^{(l,j)}\circ\alpha^{(\lambda_1,\lambda_2)}\right)(x,y)
\end{align*}
is readily checked for any local section $(x,y)$ of $\EE_A^{(\lambda_1,\lambda_2;D)}$. Hence the diagram \eqref{eq:orig-1-6} commutes. Put
\begin{align*}
A_0:=A\!\left[\dfrac{1}{\lambda_1},\dfrac{1}{\lambda_2}\right].
\end{align*}
Then, since $\alpha^{(\lambda_1^{p^l},\lambda_2^{p^l})}$ is an isomorphism over $A_0$, we have
\begin{align*}
\Psi_{\uu,\vv}^{(l,j)}=\left(\alpha^{(\lambda_1^{p^l},\lambda_2^{p^l})}\right)^{-1}\circ\Theta^{(l,j)}\circ\alpha^{(\lambda_1,\lambda_2)}.
\end{align*}
Hence $\Psi_{\uu,\vv}^{(l,j)}$ is a homomorphism of group schemes over $A_0$. Equivalently, after tensoring with $A_0$, the homomorphism
\begin{align*}
\left(\Psi_{\uu,\vv}^{(l,j)}\right)^\ast\otimes_A A_0:
\mathcal{O}\left(\EE_A^{(\lambda_1^{p^l},\lambda_2^{p^l};D')}\right)\otimes_A A_0
\longrightarrow
\mathcal{O}\left(\EE_A^{(\lambda_1,\lambda_2;D)}\right)\otimes_A A_0
\end{align*}
is a homomorphism of Hopf $A_0$-algebras, where $\mathcal{O}(\EE_A^{(\lambda_1,\lambda_2;D)})$ (resp. $\mathcal{O}(\EE_A^{(\lambda_1^{p^l},\lambda_2^{p^l};D')})$) denotes the coordinate ring of $\EE_A^{(\lambda_1,\lambda_2;D)}$ (resp. $\EE_A^{(\lambda_1^{p^l},\lambda_2^{p^l};D')}$). These coordinate rings are flat over $A$, and $A\to A_0$ is injective because $\lambda_1$ and $\lambda_2$ are nonzero divisors. Therefore the natural map
\begin{align*}
i:\mathcal{O}\left(\EE_A^{(\lambda_1,\lambda_2;D)}\right)\otimes_A\mathcal{O}\left(\EE_A^{(\lambda_1,\lambda_2;D)}\right)
\hookrightarrow
\left(\mathcal{O}\left(\EE_A^{(\lambda_1,\lambda_2;D)}\right)\otimes_A\mathcal{O}\left(\EE_A^{(\lambda_1,\lambda_2;D)}\right)\right)\left[\dfrac{1}{\lambda_1},\dfrac{1}{\lambda_2}\right]
\end{align*}
is injective.  We therefore have
\begin{align*}
i\circ\left(\left(\left(\Psi_{\uu,\vv}^{(l,j)}\right)^\ast\otimes \left(\Psi_{\uu,\vv}^{(l,j)}\right)^\ast\right)\circ \left(m_{\EE'}\right)^\ast\right)
=i\circ\left(\left(m_\EE\right)^\ast\circ\left(\Psi_{\uu,\vv}^{(l,j)}\right)^\ast\right),
\end{align*}
where $(m_\EE)^\ast$ (resp. $(m_\EE')^\ast$) is the comultiplication of $\EE_A^{(\lambda_1,\lambda_2;D)}$ (resp. $\EE_A^{(\lambda_1^{p^l},\lambda_2^{p^l};D')}$). Since $i$ is injective, this shows that $\Psi_{\uu,\vv}^{(l,j)}$ is a homomorphism of group schemes over $A$.

We give the module-theoretic argument for finiteness explicitly.  Write
$X',Y'$ for the coordinates on the target and put
\begin{align*}
R'&:=\mathcal O\bigl(\mathcal E_A^{(\lambda_1^{p^l},
\lambda_2^{p^l};D')}\bigr),\\
R_1&:=R'[X]/\bigl(\psi_1(X)-X'\bigr).
\end{align*}
Since $\psi_1(X)-X'$ is monic of degree $p^l$, $R_1$ is a free
$R'$-module with basis $1,X,\ldots,X^{p^l-1}$.  In $R_1$ we have
\begin{align*}
(1+\lambda_1X)^{p^l}=1+\lambda_1^{p^l}X'.
\end{align*}
The right-hand side is a unit in $R'$, and consequently $1+\lambda_1X$
is a unit in $R_1$.  Hence
\begin{align*}
(1+\lambda_1X)^j\bigl(\psi_2(X,Y)-Y'\bigr)
\end{align*}
is a monic polynomial of degree $p^l$ in $Y$.  It follows that
\begin{align*}
R_2:=R_1[Y]/\bigl((1+\lambda_1X)^j(\psi_2(X,Y)-Y')\bigr)
\end{align*}
is free over $R_1$ with basis $1,Y,\ldots,Y^{p^l-1}$.  Moreover, in
$R_2$ the defining equation gives
\begin{align*}
(D(X)+\lambda_2Y)^{p^l}
=(1+\lambda_1X)^j\bigl(D'(X')+\lambda_2^{p^l}Y'\bigr).
\end{align*}
Both factors on the right-hand side are units.  Therefore
$D(X)+\lambda_2Y$ is already a unit in $R_2$, so the localization occurring
in the coordinate ring of the source adds nothing.  Thus
\begin{align*}
R_2\simeq\mathcal O\bigl(\mathcal E_A^{(\lambda_1,\lambda_2;D)}\bigr)
\end{align*}
as $R'$-algebras.  The preceding two monic presentations therefore give
the basis
\begin{align*}
\left\{X^aY^b\ \middle|\ 0\leq a,b<p^l\right\}
\end{align*}
of $\mathcal{O}\bigl(\mathcal{E}_A^{(\lambda_1,\lambda_2;D)}\bigr)$ over $\mathcal{O}\bigl(\mathcal{E}_A^{(\lambda_1^{p^l},\lambda_2^{p^l};D')}\bigr)$. Consequently, $\Psi_{\uu,\vv}^{(l,j)}$ is finite locally free of rank $p^{2l}$. Since this rank is positive, it is faithfully flat and hence surjective.
\end{proof}
\begin{rem}\label{st:orig-4-12}
{\rm
Let $\zeta_2$ be a primitive $p^2$-th root of unity, and put $A:=\ZZ_{(p)}[\zeta_2]$. Note that $A$ is a discrete valuation ring with uniformizer $\lambda_2:=\zeta_2-1$. Set $\lambda:=\zeta_2^p-1$. Moreover, put
\begin{align*}
\eta:=\sum_{k=1}^{p-1}\frac{(-1)^{k-1}}{k}\lambda_2^k\quad \text{and}\quad
\widetilde{\eta}:=\frac{\lambda^{p-1}}{p}(p\eta-\lambda).
\end{align*}
Note that $[\eta]\in\widehat{W}(A_{\lambda})^{F^{(\lambda)}}$ and $[\widetilde{\eta}]\in\widehat{W}(A_{\lambda^p})^{F^{(\lambda^p)}}$. Let $\mathcal{E}_A^{(\lambda,\lambda;D)}$ be the extension associated
with $[\eta]$, and let $\mathcal{E}_A^{(\lambda^p,\lambda^p;D')}$ be the extension associated with $[\widetilde{\eta}]$. As usual, we denote $\mathcal{E}_A^{(\lambda,\lambda;D)}$ and $\mathcal{E}_A^{(\lambda^p,\lambda^p;D')}$ by $\mathcal{W}_2$ and $\mathcal{V}_2$, respectively. Put $\aa_1:=p[\lambda]/\lambda^p\in W(A)$. Since the congruence
\begin{align*}
p[\eta]-[\lambda] \equiv T_{\aa_1}([\widetilde{\eta}]) \pmod{\lambda^p}
\end{align*}
holds, the homomorphism
\begin{align*}
\Psi:=\Psi_{[\eta],[\widetilde{\eta}]}^{(1,1)}:\mathcal{W}_2\rightarrow\mathcal{V}_2
\end{align*}
is well-defined. Moreover, its kernel $\mathcal{N}_A$ is isomorphic to the constant group scheme $(\ZZ/p^2\ZZ)_A$.  Hence the short exact sequence
\begin{equation*}
\begin{tikzcd}
0 \arrow[r] &
\mathcal{N}_A\simeq(\ZZ/p^2\ZZ)_A \arrow[r] &
\mathcal{W}_2 \arrow[r,"\Psi"] &
\mathcal{V}_2 \arrow[r] &
0
\end{tikzcd}
\end{equation*}
is the Kummer--Artin--Schreier--Witt exact sequence of degree $p^2$ (see \cite[Theorem~5.2 and Example~5.3]{SS1}; see also \cite{GM,SS2}).
}
\end{rem}
\begin{rem}\label{st:orig-4-13}
{\rm
Assume that $A$ is an $\FF_p$-algebra. Let $\lambda_1$ and $\lambda_2$ be nonzero divisors of $A$, and assume that the image of $\lambda_1$ in $A_{\lambda_2}$ is nilpotent. Take $\uu\in\widehat{W}(A_{\lambda_2})^{F^{(\lambda_1)}}$. Let $\widetilde{\uu}\in \widehat{W}(A_{\lambda_2^{p^l}})$ be a lift of $\uu$. Then we have
\begin{align*}
F^{(\lambda_1)}(\widetilde{\uu})=(\alpha_0\lambda_2,\alpha_1\lambda_2,\ldots)\in \widehat{W}(A_{\lambda_2^{p^l}}).
\end{align*}
Hence $\left(F^{(\lambda_1)}(\widetilde{\uu})\right)^{(p^l)}=0$ in $\widehat{W}(A_{\lambda_2^{p^l}})$. Put $\vv:=\widetilde{\uu}^{(p^l)}\in \widehat{W}(A_{\lambda_2^{p^l}})$. Then
\begin{align*}
F^{(\lambda_1^{p^l})}(\vv)
&=\left(F-[\lambda_1^{p^l(p-1)}]\right)\left(\widetilde{\uu}^{(p^l)}\right)
=F(\widetilde{\uu})^{(p^l)}-\left([\lambda_1^{p-1}]\widetilde{\uu}\right)^{(p^l)}\\
&=\left(F(\widetilde{\uu})-[\lambda_1^{p-1}]\widetilde{\uu}\right)^{(p^l)}
=\left(F^{(\lambda_1)}(\widetilde{\uu})\right)^{(p^l)}
=0\in \widehat{W}(A_{\lambda_2^{p^l}}).
\end{align*}
Therefore, we obtain $\vv\in\widehat{W}(A_{\lambda_2^{p^l}})^{F^{(\lambda_1^{p^l})}}$. By \cite[proof of Lemma~1, pp.~123--124]{A1}, we have $T_{\aa_1}=V^l$. Note that, since $A$ is an $\FF_p$-algebra, $\widetilde{\uu}^{(p^l)}=F^l(\widetilde{\uu})$ holds. Since multiplication by $p^l$ on the Witt vectors is given by $V^lF^l$, we obtain
\begin{align*}
p^l\widetilde{\uu}
=
V^lF^l(\widetilde{\uu})
=
V^l(\vv)
=
T_{\aa_1}(\vv).
\end{align*}
Hence, if $j=0$, the congruence {\rm \eqref{eq:orig-1-2}} holds. Choose a normalized lift
\begin{align*}
D(X)=\sum_n a_nX^n\in A[X]
\end{align*}
of $E_p(\widetilde{\uu},\lambda_1;X)\in A_{\lambda_2^{p^l}}[X]$, and put
\begin{align*}
D'(X):=\sum_n a_n^{p^l}X^n\in A[X].
\end{align*}
By the functoriality of the deformed Artin--Hasse exponential with respect to homomorphisms of coefficient rings, $D'(X)$ is a lift of $E_p(\vv,\lambda_1^{p^l};X)\in A_{\lambda_2^{p^l}}[X]$. Moreover, we have
\begin{align*}
D'(X^{p^l})
=\sum_n a_n^{p^l}X^{np^l}
=\left(\sum_n a_nX^n\right)^{p^l}
=D(X)^{p^l}.
\end{align*}
Since $A$ is an $\FF_p$-algebra, we also have
\begin{align*}
\psi_1(X)
=\dfrac{1}{\lambda_1^{p^l}}\left\{(1+\lambda_1X)^{p^l}-1\right\}=X^{p^l}.
\end{align*}
Hence we obtain
\begin{align*}
D'(\psi_1(X))=D'(X^{p^l})=D(X)^{p^l}.
\end{align*}
If $j=0$, it follows that
\begin{align*}
\psi_2(X,Y)=\frac{1}{\lambda_2^{p^l}}\left\{(D(X)+\lambda_2Y)^{p^l}-D'(\psi_1(X))\right\}
=Y^{p^l}.
\end{align*}
Consequently, the equality
\begin{align*}
\Psi_{\uu,\vv}^{(l,0)}(X,Y)=\left(X^{p^l},Y^{p^l}\right)
\end{align*}
holds. Thus, with the above choices of $D$ and $D'$, $\Psi_{\uu,\vv}^{(l,0)}$ is the $l$-th Frobenius-type homomorphism.
}
\end{rem}

\subsection{The formal completion of Kummer--Artin--Schreier--Witt type isogeny}

Let $B$ be an $A$-algebra. We assume throughout this subsection that, for each $i=1,2$, every $\nu_{i,k}$ $(0\leq k\leq l-1)$ is nilpotent in $B$. We first show that, under this assumption, the kernel of $\Psi_{\uu,\vv}^{(l,j)}$ is supported on the formal neighborhood of the zero section.

\begin{lem}\label{st:orig-4-14}
The augmentation ideal of $\mathcal O(\mathcal N_B)$ is
nilpotent. In particular,
\begin{align*}
\widehat{\mathcal N}_B
=\mathcal N_B
=\Ker\left[\Psi_{\uu,\vv}^{(l,j)}:\mathcal E_B^{(\lambda_1,\lambda_2;D)}\rightarrow\mathcal E_B^{(\lambda_1^{p^l},\lambda_2^{p^l};D')}\right].
\end{align*}
Here $\mathcal O(\mathcal N_B)$ is the coordinate ring of $\mathcal N_B$.
\end{lem}

\begin{proof}
For $i=1,2$, we put
\begin{align*}
\psi_i^{(l)}(T)
:=\dfrac{1}{\lambda_i^{p^l}}\left\{(1+\lambda_iT)^{p^l}-1\right\}
=T^{p^l}+\sum_{r=1}^{p^l-1}\alpha_{i,r}T^r.
\end{align*}
By the definition of $\alpha_{i,r}$ and the nilpotency assumption on the $\nu_{i,k}$, each $\alpha_{i,r}$ is nilpotent in $B$. Let $\overline{X}$ and $\overline{Y}$ denote the images of $X$
and $Y$ in $\mathcal O(\mathcal N_B)$. Set
\begin{align*}
I_1:=
(\alpha_{1,1},\ldots,\alpha_{1,p^l-1})B.
\end{align*}
Since $I_1$ is generated by finitely many nilpotent elements, $I_1$ is a nilpotent ideal. By
\begin{align*}
\psi_1(\overline{X})=\overline{X}^{p^l}+\sum_{r=1}^{{p^l}-1}\alpha_{1,r}\overline{X}^r=0,
\end{align*}
we obtain $\overline{X}^{p^l}\in I_1\mathcal O(\mathcal N_B)$. Hence $\overline{X}$ is nilpotent. Next, put
\begin{align}\label{eq:orig-4-6}
h(X)
:=\psi_2(X,0)=\frac{1}{\lambda_2^{p^l}}\left\{\frac{D(X)^{p^l}}{(1+\lambda_1X)^j}-D'(\psi_1(X))\right\}.
\end{align}
Since $\Psi_{\uu,\vv}^{(l,j)}$ preserves the zero section, we have $h(0)=0$. Therefore $h(\overline{X})$ is nilpotent. Then we have
\begin{align}\label{eq:orig-4-7}
\begin{split}
\psi_2(X,Y)-h(X)
&=\dfrac{1}{\lambda_2^{p^l}}\left\{\dfrac{(D(X)+\lambda_2Y)^{p^l}-D(X)^{p^l}}{(1+\lambda_1 X)^j}\right\}\\
&=\dfrac{1}{(1+\lambda_1X)^j}
\left\{Y^{p^l}+\sum_{r=1}^{p^l-1}
\alpha_{2,r}D(X)^{p^l-r}Y^r\right\}.
\end{split}
\end{align}
Since $\psi_2(\overline X,\overline Y)=0$ on $\mathcal N_B$, \eqref{eq:orig-4-7} gives
\begin{align*}
\overline{Y}^{p^l}
+\sum_{r=1}^{p^l-1}
\alpha_{2,r}D(\overline X)^{p^l-r}\overline Y^r
=-(1+\lambda_1\overline X)^j h(\overline X).
\end{align*}
Put
\begin{align*}
I_2
:=
(h(\overline X),\alpha_{2,1},\ldots,\alpha_{2,p^l-1})
\mathcal O(\mathcal N_B).
\end{align*}
Then $I_2$ is a nilpotent ideal. Hence we obtain $\overline{Y}^{p^l}\in I_2$. Therefore, $\overline{Y}$ is nilpotent. Thus the augmentation ideal $(\overline{X},\overline{Y})\subset\mathcal O(\mathcal N_B)$ is nilpotent. Consequently, the completion of $\mathcal N_B$ along the zero section coincides with $\mathcal N_B$.
\end{proof}
We shall also use the following immediate consequence.  Since the augmentation
ideal of $\mathcal N_B$ is nilpotent, every homomorphism
$\mathcal N_B\to\mathbb G_{m,B}$ factors uniquely through the formal completion
$\widehat{\mathbb G}_{m,B}$.  Hence, after every base change $B\to C$, we have
\begin{align}\label{eq:orig-4-8}
\Hom(\mathcal N_C,\widehat{\mathbb G}_{m,C})
=\Hom(\mathcal N_C,\mathbb G_{m,C}).
\end{align}
We now consider the formal completion of $\Psi_{\uu,\vv}^{(l,j)}$ along the zero sections.

\begin{lem}\label{st:orig-4-15}
Under the above assumptions, the induced homomorphism
\begin{align*}
\widehat{\Psi}_{\uu,\vv}^{(l,j)}:
\widehat{\mathcal E}_B^{(\lambda_1,\lambda_2;D)}
\rightarrow
\widehat{\mathcal E}_B^{(\lambda_1^{p^l},\lambda_2^{p^l};D')}
\end{align*}
is finite faithfully flat. In particular, it is an epimorphism for the fppf topology. Moreover, the induced pull-back
\begin{align*}
\left(\widehat{\Psi}_{\uu,\vv}^{(l,j)}\right)^\ast
:\mathcal{O}\left(\widehat{\mathcal E}_B^{(\lambda_1^{p^l},\lambda_2^{p^l};D')}\right)
\rightarrow
\mathcal{O}\left(\widehat{\mathcal E}_B^{(\lambda_1,\lambda_2;D)}\right)
\end{align*}
is injective.
\end{lem}

\begin{proof}
Put
\begin{align*}
R
:=\mathcal{O}\left(\mathcal E_B^{(\lambda_1,\lambda_2;D)}\right)\quad\text{and}\quad
R'
:=\mathcal{O}\left(\mathcal E_B^{(\lambda_1^{p^l},\lambda_2^{p^l};D')}\right).
\end{align*}
Let $I\subset R$ and $I'\subset R'$ be the augmentation ideals. Since $\Psi_{\uu,\vv}^{(l,j)}$ preserves the zero sections, we have $I'R\subset I$. Moreover, we have
\begin{align*}
R/I'R\simeq \mathcal{O}(\mathcal N_B).
\end{align*}
By Lemma~\ref{st:orig-4-14}, the augmentation ideal $I/I'R$ of $\mathcal O(\mathcal N_B)$ is nilpotent. Hence there exists an integer $N>0$ such that $I^N\subset I'R$. Thus the $I$-adic topology and the $I'R$-adic topology on $R$ are equivalent. Since $\Psi_{\uu,\vv}^{(l,j)}$ is finite locally free, $R$ is a finite locally free $R'$-module. Since completion commutes with finite locally free modules, we obtain
\begin{align*}
\widehat{R}_I
\simeq
\widehat{R}_{I'R}
\simeq
R\otimes_{R'}\widehat{R'}_{I'},
\end{align*}
where $\widehat{R}_I=\varprojlim_{n} R/I^n$, $\widehat{R}_{I'R}=\varprojlim_{n} R/(I'R)^n$ and $\widehat{R'}_{I'}=\varprojlim_{n} R'/{I'}^n$. Consequently, the homomorphism $\widehat{R'}_{I'}\rightarrow\widehat{R}_I$ is obtained from the faithfully flat homomorphism $R'\rightarrow R$ by base change. Hence it is finite faithfully flat and, in particular, injective. It follows that $\widehat{\Psi}_{\uu,\vv}^{(l,j)}$ is finite faithfully flat, and hence an epimorphism for the fppf topology. The injectivity of the above homomorphism of completed coordinate rings gives the injectivity of the pull-back $\left(\widehat{\Psi}_{\uu,\vv}^{(l,j)}\right)^\ast$.
\end{proof}
As a consequence, we obtain an exact sequence of fppf sheaves
\begin{align*}
0
\rightarrow
\widehat{\mathcal N}_B
\rightarrow
\widehat{\mathcal E}_B^{(\lambda_1,\lambda_2;D)}
\xrightarrow{\widehat{\Psi}_{\uu,\vv}^{(l,j)}}
\widehat{\mathcal E}_B^{(\lambda_1^{p^l},\lambda_2^{p^l};D')}
\rightarrow
0.
\end{align*}
Since $\widehat{\mathcal N}_B=\mathcal N_B$, this may also be written as
\begin{align*}
0
\rightarrow
\mathcal N_B
\rightarrow
\widehat{\mathcal E}_B^{(\lambda_1,\lambda_2;D)}
\xrightarrow{\widehat{\Psi}_{\uu,\vv}^{(l,j)}}
\widehat{\mathcal E}_B^{(\lambda_1^{p^l},\lambda_2^{p^l};D')}
\rightarrow
0.
\end{align*}
\section{Cartier duals of Kummer--Artin--Schreier--Witt type kernels}
Let $A$ be a $\ZZ_{(p)}$-algebra. Take nonzero divisors $\lambda_1,\lambda_2\in A$ such that $\lambda_1^{p^l-p^{l-1}}\mid p$ and $\lambda_2^{p^l-p^{l-1}}\mid p$. Let $\EE_A^{(\lambda_1,\lambda_2;D)}$ be an extension associated with $\uu\in \widehat{W}(A_{\lambda_2})^{F^{(\lambda_1)}}$, and let $\widehat{\EE}_A^{(\lambda_1,\lambda_2;D)}$ denote the formal completion of $\EE_A^{(\lambda_1,\lambda_2;D)}$ along the zero section. In the formal case, after transporting the extension through the canonical change-of-lift isomorphism, we may take
\begin{align*}
D=D(X)=E_p(\uu^\#,\lambda_1;X)\in A[[X]],
\end{align*}
where $\uu^\#\in W(A)$ is the lift chosen in the Introduction (for the
general case, see \cite[Section~4]{MRT2}).  More explicitly, if $D_0$ is
the normalized polynomial lift used in Section~4, write
\begin{align*}
D_0(X)-D(X)=\lambda_2h(X),\qquad h(0)=0.
\end{align*}
Then the change of coordinates
\begin{align*}
\theta_h:\widehat{\mathcal E}^{(\lambda_1,\lambda_2;D_0)}_A
\xrightarrow{\ \sim\ }
\widehat{\mathcal E}^{(\lambda_1,\lambda_2;D)}_A,\qquad
(X,Y)\longmapsto (X,Y+h(X)),
\end{align*}
preserves the multiplicative coordinates
$(1+\lambda_1X,D_0(X)+\lambda_2Y)$.  Likewise, if
$D'_0-D'=\lambda_2^{p^l}h'$, we use the analogous isomorphism
$\theta_{h'}$ for the target.  The transported formal homomorphism satisfies
\begin{align*}
\theta_{h'}\circ\widehat\Psi_{D_0,D'_0}
=\widehat\Psi_{D,D'}\circ\theta_h,
\end{align*}
where the subscripts record the chosen lifts, because both sides induce the
same homomorphism on the two multiplicative
coordinates.  Indeed, since $\lambda_1$ and $\lambda_2$ are nonzero
divisors in $A$, these two multiplicative coordinates determine the two
formal coordinates.  Thus the kernel and its Cartier dual are unchanged
by this replacement.
\begin{rem}\label{st:orig-5-1}
{\rm
Throughout this paper, we consider only those extension classes in
\begin{align*}
\Ext^1\left(\widehat{\mathcal{G}}_A^{(\lambda_1)},\widehat{\mathcal{G}}_A^{(\lambda_2)}\right)
\end{align*}
that arise from homomorphisms $\widehat{\mathcal{G}}^{(\lambda_1)}_{A_{\lambda_2}}\rightarrow\widehat{\mathbb{G}}_{m,A_{\lambda_2}}$. Since the isomorphism class of $\widehat{\EE}_A^{(\lambda_1,\lambda_2;D)}$ depends only on $
\uu\in\widehat{W}(A_{\lambda_2})^{F^{(\lambda_1)}}$, if we choose a polynomial lift $\widetilde{D}\in A[X]$ of
$E_p(\uu,\lambda_1;X)$, then there exists an isomorphism
\begin{align*}
\widehat{\EE}_A^{(\lambda_1,\lambda_2;D)}
\simeq
\widehat{\EE}_A^{(\lambda_1,\lambda_2;\widetilde{D})}.
\end{align*}
}
\end{rem}

\subsection{Key Lemmas}

We use throughout the fixed lifts $\uu^\#$ and $\widetilde{\vv}$ from
the Introduction. Put
\begin{align*}
D(X)&=E_p(\uu^\#,\lambda_1;X),\\
D'(X)&=E_p(\widetilde{\vv},\lambda_1^{p^l};X),\\
\dd&=F^{(\lambda_1^{p^l})}(\widetilde{\vv})
  =(\lambda_2^{p^l}b'_0,\lambda_2^{p^l}b'_1,\ldots),\\
E_0(X)&=E_p(p^l[\lambda_1],\lambda_1;X)=(1+\lambda_1X)^{p^l}.
\end{align*}
The series $G_p(\boldsymbol\eta,\lambda_1^{p^l};E_0)$ is the
integral one-dimensional correction of Lemma~\ref{st:orig-4-4}.
For $r\geq0$, define
\begin{align*}
E_{1,r}(X)
&:=E_p\!\left(
\widetilde p^{\,r}(p^l\uu^\#-T_{\aa_1}\widetilde{\vv}-j[\lambda_1]),
\lambda_1;X\right),\\
E_{2,r}(X,Y)
&:=E_p\!\left(\widetilde p^{\,r}(p^l[\lambda_2]),
\lambda_2;\frac{Y}{D(X)}\right),\\
K_r(X,Y)
&:=\frac{E_{1,r}(X)E_{2,r}(X,Y)}
{G_p(\widetilde p^{\,r}\dd,\lambda_1^{p^l};E_0)}.
\end{align*}
Each displayed exponential retains its indicated Witt-vector
presentation. Write $E_1=E_{1,0}$ and $E_2=E_{2,0}$.
Lemma~\ref{st:orig-4-4} gives
\begin{align*}
D'(\psi_1(X))
=E_p(T_{\aa_1}\widetilde{\vv},\lambda_1;X)
G_p(\dd,\lambda_1^{p^l};E_0).
\end{align*}
Using the formula for the transported $\psi_2$, we therefore obtain
\begin{align}\label{eq:orig-5-1}
K_0(X,Y)=\frac{E_1E_2}{G_p(\dd,\lambda_1^{p^l};E_0)}
=1+\lambda_2^{p^l}\frac{\psi_2(X,Y)}{D'(\psi_1(X))}.
\end{align}
This identity also specifies the integral series denoted by
$(K_0-1)/\lambda_2^{p^l}$: it is $\psi_2/D'(\psi_1)$ before any
base change.

With these Witt-vector presentations, define the second correction by
\begin{align}\label{eq:second-correction}
G_p(\boldsymbol\eta,\lambda_2^{p^l};K_0)
:=\prod_{r\geq1}
\left(
\frac{1+(K_0(X,Y)-1)^{p^r}}{K_r(X,Y)}
\right)^{\Phi_{r-1}(\boldsymbol\eta)/(p^r\lambda_2^{p^{l+r}})}.
\end{align}
The rational exponents in this display are interpreted by specialization
of the universal integral coefficients, as in
\eqref{eq:one-dimensional-correction}. The specialization for the present
data is explained in Lemma~\ref{st:orig-5-2}.
The product is well-defined in formal degree. In fact,
$E_{i,r}$ is $1$ modulo $(X,Y)^{p^r}$, and
$G_p(\widetilde p^{\,r}\dd,\lambda_1^{p^l};E_0)$ is $1$ modulo
$X^{p^{r+1}}$, so $K_r$ is $1$ modulo $(X,Y)^{p^r}$.
For $\ww=(\ww_1,\ww_2)$ put
\begin{align}\label{eq:gamma-definition}
\Gamma_{\ww}(X,Y)
:=G_p(\ww_1,\lambda_1^{p^l};E_0)
G_p(\ww_2,\lambda_2^{p^l};K_0).
\end{align}
In \eqref{eq:second-correction} the entire sequence $(K_r)_{r\geq0}$
is part of the definition. In particular, a substitution taking
$K_0$ to $1$ need not take all the $K_r$ to $1$.

We use the method of \cite[Propositions~4.17 and~9.3]{SS2} with two deformation parameters.
The pullback formula is analogous to \cite[(98), p.~75]{SS2}.

\begin{lem}\label{st:orig-5-2}
The series $G_p(\boldsymbol\eta,\lambda_2^{p^l};K_0)$ has coefficients
in $A[\boldsymbol\eta]$, defined by specialization of universal integral
polynomials. Consequently, for every $A$-algebra
$C$, \eqref{eq:gamma-definition} defines
\begin{align*}
\Gamma_{\ww}\in1+(X,Y)C[[X,Y]]\qquad(\ww\in W^2(C)).
\end{align*}
These series commute with base change and satisfy
\begin{align*}
\Gamma_{\ww+\ww'}=\Gamma_{\ww}\Gamma_{\ww'},
\qquad \Gamma_0=1.
\end{align*}
For every $\xx\in W^2(C)$, we have
\begin{align}\label{eq:corrected-pullback}
\begin{split}
&E_p(\xx,(\lambda_1^{p^l},\lambda_2^{p^l});\Psi(X,Y))\\
&\qquad=
E_p(\Upsilon_2\xx,(\lambda_1,\lambda_2);X,Y)
\Gamma_{U_2'\xx}(X,Y).
\end{split}
\end{align}
\end{lem}

\begin{proof}
Let $\Lambda_1,\Lambda_2,M,N$ and the coordinates of the Witt vectors
$\mathbf A,\mathbf B,\mathbf C,\mathbf E$ be independent variables over
$\ZZ_{(p)}$, and let $S$ be the polynomial ring in these variables.
Here $M\mathbf A=(MA_0,MA_1,\ldots)$, and similarly for the other
vectors. For formal variables $X,Z$, put
\begin{align*}
\mathcal E_0&=E_p(M\mathbf A,\Lambda_1;X),\\
\mathcal K_r&=
\frac{E_p(\widetilde p^{\,r}(N\mathbf C),\Lambda_1;X)
      E_p(\widetilde p^{\,r}(N\mathbf E),\Lambda_2;Z)}
     {G_p(\widetilde p^{\,r}(N\mathbf B),M;\mathcal E_0)}
\qquad(r\geq0).
\end{align*}
By \cite[Proposition~4.12]{SS2}, the inner $G_p$ has coefficients
in $S$, since the componentwise quotient of $M\mathbf A$ by $M$ is
$\mathbf A$. It has constant term $1$. Thus $\mathcal K_r\in S[[X,Z]]$.
Since $\mathcal K_0\equiv1\bmod N$, we have
$T:=(\mathcal K_0-1)/N\in(X,Z)S[[X,Z]]$.

For a Witt vector of variables $\boldsymbol\eta$, put
$P=S[\boldsymbol\eta]$. Let $G_p(\boldsymbol\eta,N;\mathcal K_0)$ denote the product in
\eqref{eq:second-correction} with parameter $N$ and sequence
$(\mathcal K_r)_{r\geq0}$, first taken over $\operatorname{Frac}(P)$.
Only finitely many factors contribute to each coefficient, since
$\mathcal K_r-1\in(X,Z)^{p^r}S[[X,Z]]$.
For another Witt vector of variables $\zz$, put $P'=S[\zz]$ and
\begin{align*}
\beta_0=\frac{z_0}{N},\qquad
\beta_r=\frac{\Phi_{r-1}(F^{(N)}\zz)}{p^rN^{p^r}}\quad(r\geq1).
\end{align*}
Applying \cite[Propositions~4.11 and~4.13]{SS2} to the two
exponentials and to the inner $G_p$, respectively, we obtain over
$\operatorname{Frac}(P')$
\begin{align*}
\prod_{r\geq0}\mathcal K_r^{\beta_r}
&=\frac{E_p(T_{\mathbf C}\zz,\Lambda_1;X)
        E_p(T_{\mathbf E}\zz,\Lambda_2;Z)}
       {G_p(T_{\mathbf B}\zz,M;\mathcal E_0)}.
\end{align*}
Here $M$ and $N$ are the two parameters in Proposition~4.13.
By the product formula for $E_p$ and $1+NT=\mathcal K_0$,
\begin{align*}
G_p(F^{(N)}\zz,N;\mathcal K_0)
&=\frac{E_p(\zz,N;T)}{\prod_{r\geq0}\mathcal K_r^{\beta_r}}\\
&=\frac{E_p(\zz,N;T)\,G_p(T_{\mathbf B}\zz,M;\mathcal E_0)}
       {E_p(T_{\mathbf C}\zz,\Lambda_1;X)
        E_p(T_{\mathbf E}\zz,\Lambda_2;Z)}
\in P'[[X,Z]].
\end{align*}
Indeed, the two series in the denominator are units in $P'[[X,Z]]$.
By \cite[Corollary~1.8]{SS1}, the homomorphism
$P\to P'$, $\boldsymbol\eta\mapsto F^{(N)}\zz$, is faithfully flat.
If $a/d\in\operatorname{Frac}(P)$ is a coefficient of
$G_p(\boldsymbol\eta,N;\mathcal K_0)$, the last equality gives
$a\in dP'\cap P=dP$. Hence $G_p(\boldsymbol\eta,N;\mathcal K_0)\in P[[X,Z]]$.

We now specialize
\begin{align*}
(\Lambda_1,\Lambda_2,M,N)&\longmapsto
(\lambda_1,\lambda_2,\lambda_1^{p^l},\lambda_2^{p^l}),\\
(\mathbf A,\mathbf B,\mathbf C,\mathbf E)&\longmapsto
(\aa_1,\bb',\cc,\aa_2),\qquad Z\longmapsto Y/D(X).
\end{align*}
Then $\mathcal E_0,\mathcal K_r$ specialize to $E_0,K_r$.
By \eqref{eq:orig-5-1} and the fact that $\lambda_2$ is a nonzero
divisor in $A$, $T$ specializes to $\psi_2/D'(\psi_1)$.
This proves the integrality assertion. For both correction factors,
additivity of the ghost components gives multiplicativity in the Witt-vector
argument over the fraction field, hence over the integral coefficient ring.
It follows by specialization that
$\Gamma_{\ww+\ww'}=\Gamma_{\ww}\Gamma_{\ww'}$ and $\Gamma_0=1$,
also after any base change $A\to C$.
The above equality for $G_p(F^{(N)}\zz,N;\mathcal K_0)$ gives
\begin{align}\label{eq:second-correction-descent}
\begin{split}
&G_p(F^{(\lambda_2^{p^l})}\zz,\lambda_2^{p^l};K_0)\\
&\quad=
\frac{
E_p(\zz,\lambda_2^{p^l};\psi_2/D'(\psi_1))
G_p(T_{\bb'}\zz,\lambda_1^{p^l};E_0)
}{
E_p(T_{\cc}\zz,\lambda_1;X)
E_p(T_{\aa_2}\zz,\lambda_2;Y/D(X))
}.
\end{split}
\end{align}
Finally, rearrange \eqref{eq:second-correction-descent}, take
$\zz=\xx_2$, and multiply by the identity of
Lemma~\ref{st:orig-4-4} for $\xx_1$ and $\lambda_1$.
The ratio of the one-dimensional correction factors is
\begin{align*}
\frac{G_p(F^{(\lambda_1^{p^l})}\xx_1,\lambda_1^{p^l};E_0)}
{G_p(T_{\bb'}\xx_2,\lambda_1^{p^l};E_0)}
=G_p(F^{(\lambda_1^{p^l})}\xx_1-T_{\bb'}\xx_2,
\lambda_1^{p^l};E_0).
\end{align*}
The two arguments of the resulting correction factors are exactly
the components of $U_2'\xx$. This proves
\eqref{eq:corrected-pullback} as an identity of specialized universal
integral series, and hence over every $A$-algebra $C$.
\end{proof}

\begin{sublem}\label{st:orig-5-6}
The morphisms $U_2,U_2':W_A^2\to W_A^2$ are faithfully flat.
In particular they are epimorphisms for the fpqc topology.
\end{sublem}
\begin{proof}
Each $F^{(\lambda)}:W_A\to W_A$ is faithfully flat by
\cite[Corollary~1.8]{SS1}.
The morphism $U_2$ is the composite of $F^{(\lambda_1)}\times
\operatorname{id}$, the automorphism
$(\xx,\yy)\mapsto(\xx-T_{\bb}\yy,\yy)$, and
$\operatorname{id}\times F^{(\lambda_2)}$.
Thus it is faithfully flat. The same argument applies to $U_2'$.
These morphisms are affine, hence quasi-compact, which gives the
fpqc assertion.
\end{proof}

Hereafter let $B$ be an $A$-algebra in which every $\nu_{i,k}$
is nilpotent.

\begin{lem}\label{st:orig-5-7}
There is a unique group scheme homomorphism
$\Upsilon_2':W_A^2\to W_A^2$ satisfying
\begin{align}\label{eq:upsilon-factorization}
U_2\Upsilon_2=\Upsilon_2'U_2'.
\end{align}
Its base change to $B$ is a monomorphism.
\end{lem}
\begin{proof}
Let $C$ be any $A$-algebra and take $\xx\in\Ker(U_2')(C)$.
By \eqref{eq:orig-3-2}, $E_p(\xx,(\lambda_1^{p^l},\lambda_2^{p^l}))$ is a character
of the target formal group. Since $\Gamma_0=1$,
\eqref{eq:corrected-pullback} identifies its pullback with
$E_p(\Upsilon_2\xx,(\lambda_1,\lambda_2))$.
By \eqref{eq:orig-3-2} and the injectivity of the two-dimensional
exponential map proved in Subsection~\ref{subsec:two-dimensional-exponentials}, this implies
$U_2\Upsilon_2\xx=0$.
Thus $U_2\Upsilon_2$ vanishes on $\Ker U_2'$.
By Sublemma~\ref{st:orig-5-6}, the quotient of $W_A^2$ by this
kernel is $W_A^2$ as an fpqc sheaf. Its universal property gives
$\Upsilon_2'$, and Yoneda gives representability of the induced
morphism. This also proves uniqueness.

Write $S_i=T'_{\aa_i}$ for the homomorphisms constructed in
Lemma~\ref{st:orig-4-5}.
The second component of \eqref{eq:upsilon-factorization}, together
with $F^{(\lambda_2)}T_{\aa_2}=S_2F^{(\lambda_2^{p^l})}$, shows that the
second component of $\Upsilon_2'(\ww_1,\ww_2)$ is $S_2\ww_2$:
one can check this after the fpqc epimorphism $U_2'$.
Applying the first component to $(\xx_1,0)$ and canceling the fpqc
epimorphism $F^{(\lambda_1^{p^l})}$ similarly gives
$\Upsilon_2'(\ww_1,0)=(S_1\ww_1,0)$.
By additivity there is a homomorphism $L:W_A\to W_A$ such that
\begin{align*}
\Upsilon_2'(\ww_1,\ww_2)
=(S_1\ww_1+L(\ww_2),S_2\ww_2).
\end{align*}
After base change to $B$, both $S_i$ are monomorphisms by
Lemma~\ref{st:orig-4-5}. If $\Upsilon_2'(\ww_1,\ww_2)=0$
over any $B$-algebra, its second component gives $\ww_2=0$ and
then its first component gives $\ww_1=0$. This proves the assertion.

We have, in particular, the commutative diagram of fpqc sheaves
\begin{equation}\label{eq:orig-5-5}
\begin{tikzcd}
0 \arrow[r] &
\Ker(U_2') \arrow[r] \arrow[d,"\Upsilon_2"'] &
W_B^2 \arrow[r,"U_2'"] \arrow[d,"\Upsilon_2"'] &
W_B^2 \arrow[r] \arrow[d,"\Upsilon_2'"] &
0\\
0 \arrow[r] &
\Ker(U_2) \arrow[r] &
W_B^2 \arrow[r,"U_2"] &
W_B^2 \arrow[r] &
0.
\end{tikzcd}
\end{equation}
\end{proof}

\begin{lem}\label{st:orig-5-3}
For every $B$-algebra $C$,
\begin{align*}
\Ker(U_2')(C)=\Ker(U_2\Upsilon_2)(C).
\end{align*}
\end{lem}
\begin{proof}
Use \eqref{eq:upsilon-factorization} and the injectivity of
$\Upsilon_2'(C)$ proved in Lemma~\ref{st:orig-5-7}.
\end{proof}

\begin{lem}\label{st:orig-5-4}
For every $B$-algebra $C$, the diagram
\begin{equation}\label{eq:orig-5-4}
\begin{tikzcd}
\Hom\!\left(\widehat{\EE}_C^{(\lambda_1^{p^l},\lambda_2^{p^l};D')},
\widehat{\GG}_{m,C}\right)
\arrow[r,"\Psi^*"] &
\Hom\!\left(\widehat{\EE}_C^{(\lambda_1,\lambda_2;D)},
\widehat{\GG}_{m,C}\right)\\
\Ker(U_2')(C)
\arrow[r,"\Upsilon_2"] \arrow[u,"\phi_1"] &
\Ker(U_2)(C) \arrow[u,"\phi_2"]
\end{tikzcd}
\end{equation}
commutes, where the vertical maps are the isomorphisms in
\eqref{eq:orig-3-2}.
\end{lem}
\begin{proof}
For $\xx\in\Ker U_2'$, \eqref{eq:corrected-pullback} has correction
$\Gamma_0=1$ and is exactly the required equality.
\end{proof}

\begin{lem}\label{st:orig-5-5}
The homomorphism $\Upsilon_2:W_B^2\to W_B^2$ is a monomorphism.
\end{lem}
\begin{proof}
Its diagonal entries $T_{\aa_1}$ and $T_{\aa_2}$ are closed
immersions by Lemma~\ref{st:orig-4-5}.
Let $C$ be a $B$-algebra and suppose that $\Upsilon_2(\xx)=0$.
The second component gives $\xx_2=0$, and the first then gives
$\xx_1=0$.
\end{proof}

\begin{rem}\label{st:orig-5-8}
{\rm
Assume that $A$ is an $\FF_p$-algebra, and keep the notation $\vv=\widetilde{\uu}^{(p^l)}=F^l(\widetilde{\uu})$ as in Remark~\ref{st:orig-4-13}. Choose a lift $\uu^{\#}\in W(A)$ of $\widetilde{\uu}\in\widehat{W}(A_{\lambda_2^{p^l}})$, and put
\begin{align*}
\widetilde{\vv}:=F^l(\uu^{\#})\in W(A).
\end{align*}
Then $\widetilde{\vv}$ is a lift of $\vv\in\widehat{W}(A_{\lambda_2^{p^l}})$. Indeed, since $A$ is an $\FF_p$-algebra, we have
\begin{align*}
\widetilde{\vv}=F^l(\uu^{\#})=\left(\uu^{\#}\right)^{(p^l)}\equiv\widetilde{\uu}^{(p^l)}=\vv\pmod{\lambda_2^{p^l}}.
\end{align*}
By \cite[proof of Lemma~1, pp.~123--124]{A1}, we have $T_{\aa_1}=T_{\aa_2}=V^l$. Moreover, since multiplication by $p^l$ on $W(A)$ is given by $V^lF^l$, we obtain
\begin{align*}
p^l\uu^{\#}=V^lF^l(\uu^{\#})=V^l(\widetilde{\vv})=T_{\aa_1}(\widetilde{\vv})\in W(A).
\end{align*}
Hence, if $j=0$, then
\begin{align*}
\cc=\frac{1}{\lambda_2^{p^l}}\left\{p^l\uu^{\#}-j[\lambda_1]-T_{\aa_1}(\widetilde{\vv})\right\}=0
\end{align*}
holds. Therefore, we have
\begin{align*}
\Upsilon_2=
\begin{pmatrix}
V^l & 0\\
0   & V^l
\end{pmatrix}
=:
V^l\times V^l.
\end{align*}
In characteristic $p$, the relation $\widetilde{\vv}=F^l(\uu^\#)$ gives $\bb'=F^l(\bb)$: apply $F^l$ to $F^{(\lambda_1)}(\uu^\#)=\lambda_2\bb$ and use componentwise Frobenius and the fact that $\lambda_2$ is a nonzero divisor in $A$. The identity $T_{\bb}V^l=V^lT_{F^l(\bb)}$ now shows that the diagram
\begin{equation*}
\begin{tikzcd}[column sep=12ex]
W_A^2
\arrow[r,"\Upsilon_2=V^l\times V^l"]
\arrow[d,"U_2'"']
&
W_A^2
\arrow[d,"U_2"]
\\
W_A^2
\arrow[r,"V^l\times V^l"]
&
W_A^2
\end{tikzcd}
\end{equation*}
is commutative. Therefore, by the uniqueness of $\Upsilon_2'$, we obtain
\begin{align*}
\Upsilon_2'=V^l\times V^l.
\end{align*}

For these choices, the Verschiebung identity for the exponential and
$D'(X^{p^l})=D(X)^{p^l}$ give, for every $A$-algebra $C$ and every
$\xx\in W^2(C)$,
\begin{align*}
E_p(\xx,(\lambda_1^{p^l},\lambda_2^{p^l});\Psi(X,Y))
=E_p((V^l\times V^l)\xx,(\lambda_1,\lambda_2);X,Y).
\end{align*}
Comparing with \eqref{eq:corrected-pullback} and using the faithful
flatness of $U_2'$, we obtain $\Gamma_{\ww}=1$ for every $\ww$.
Thus the corrected correspondence reduces to the usual exponential
correspondence in this characteristic-$p$ case.
}
\end{rem}
\subsection{Cohomology of Witt groups}\label{subsec:witt-cohomology}

We state the vanishing result for an arbitrary finite number of Witt
vectors. We shall use the case $r=2$ below.

\begin{thm}\label{st:orig-5-9}
Let $R$ be a commutative ring, put $X=\Spec R$, and let $r\geq1$.
Write
\begin{align*}
W_{n,R}^r:=\prod_{a=1}^r W_{n,R},\qquad
W_R^r:=\prod_{a=1}^r W_R=\varprojlim_n W_{n,R}^r.
\end{align*}
Here $n$ denotes the length of a $p$-typical Witt vector.
Then
\begin{align*}
H^i(X_{\mathrm{fppf}},W_{n,R}^r)=0,\qquad
H^i(X_{\mathrm{fppf}},W_R^r)=0
\end{align*}
for every $n\geq1$ and $i>0$.
\end{thm}
\begin{proof}
The fppf cohomology of a quasi-coherent module on an affine scheme
vanishes in positive degrees. Hence
\begin{align*}
H^i(X_{\mathrm{fppf}},\GG_{a,R}^r)=0\qquad(i>0).
\end{align*}
For $n\geq2$, we have an exact sequence
\begin{align*}
0\longrightarrow\GG_{a,R}^r
\xrightarrow{(V^{n-1})^{\times r}}W_{n,R}^r
\xrightarrow{R_{n-1}^{\times r}}W_{n-1,R}^r
\longrightarrow0.
\end{align*}
Since $W_{1,R}^r=\GG_{a,R}^r$, induction on $n$ gives the first
assertion. The same proof applies after replacing $R$ by any
$R$-algebra $D$.

Put $\mathcal A=\mathrm{Ab}(X_{\mathrm{fppf}})$ and
$\mathcal F_n=W_{n,R}^r$. For an affine object $U=\Spec D$ of
$X_{\mathrm{fppf}}$, we have
\begin{align*}
H^i(U_{\mathrm{fppf}},\mathcal F_n|_U)&=0\quad(i>0),\\
\mathcal F_n(U)&=W_n^r(D),\\
R_n^{\times r}:W_{n+1}^r(D)&\twoheadrightarrow W_n^r(D).
\end{align*}
The last map is truncation on each Witt vector. Thus
\begin{align*}
\varprojlim\nolimits^1_n W_n^r(D)=0.
\end{align*}

The category $\mathcal A$ has enough injectives, and products are
exact in $\mathrm{Ab}$. Apply \cite[Proposition~1.6]{Jan} to
$\Gamma(U,-):\mathcal A\to\mathrm{Ab}$. Put
$K=R\varprojlim_n\mathcal F_n$. Since $\varprojlim$ preserves
injectives \cite[Remark~1.17(a)]{Jan} and $\Gamma(U,-)$ commutes
with inverse limits, we have
\begin{align*}
R^i\!\left(\varprojlim_n\Gamma(U,-)\right)(\mathcal F_\bullet)
\simeq H^i(U_{\mathrm{fppf}},K).
\end{align*}
Consequently, for $i\geq1$, there is an exact sequence
\begin{align*}
0\longrightarrow
\varprojlim\nolimits^1_n H^{i-1}(U_{\mathrm{fppf}},\mathcal F_n|_U)
\longrightarrow H^i(U_{\mathrm{fppf}},K)
\longrightarrow
\varprojlim_n H^i(U_{\mathrm{fppf}},\mathcal F_n|_U)
\longrightarrow0.
\end{align*}
Both outer terms vanish: for $i=1$ use
$\varprojlim^1_nW_n^r(D)=0$, and for $i\geq2$ use the finite-level
vanishing. In degree zero, the definition of the derived functor gives
\begin{align*}
H^i(U_{\mathrm{fppf}},K)
=\begin{cases}
W^r(D),&i=0,\\
0,&i>0.
\end{cases}
\end{align*}

We verify that $K$ is represented by $W_R^r$ in degree zero.
Take an injective resolution $\mathcal F_\bullet\to\mathcal I_\bullet^\bullet$
in $\mathcal A^{\mathbb N}$. Then
$K^\bullet=\varprojlim_n\mathcal I_n^\bullet$ is a bounded below
complex of injectives representing $K$. Its cohomology sheaves are
the sheaves associated with
\begin{align*}
U\longmapsto H^i\!\left(\Gamma(U,K^\bullet)\right)
=H^i(U_{\mathrm{fppf}},K).
\end{align*}
Affine objects form a basis. Therefore
\begin{align*}
\mathcal H^0(K)=W_R^r,\qquad
\mathcal H^i(K)=0\quad(i>0),\qquad
W_R^r\xrightarrow{\sim}R\varprojlim_nW_{n,R}^r.
\end{align*}
Taking $U=X$ proves the second assertion.
\end{proof}

\begin{lem}\label{st:orig-5-10}
The presheaves $W_B^2/\Upsilon_2$ and $W_B^2/\Upsilon_2'$
defined in {\rm \eqref{eq:orig-1-8}} and {\rm \eqref{eq:orig-1-9}}, respectively, are fppf sheaves
on the category of $B$-algebras.
\end{lem}
\begin{proof}
Let $C$ be a $B$-algebra and put $X=\Spec\,C$. Let $\widetilde{W_B^2/\Upsilon_2}$ denote the fppf sheafification of the presheaf $W_B^2/\Upsilon_2$ on the category of $B$-algebras. Since $\Upsilon_2$ is a monomorphism, restricting
$\widetilde{W_B^2/\Upsilon_2}$ to $X_{\mathrm{fppf}}$, we have a short exact sequence of fppf sheaves
\begin{equation*}
\begin{tikzcd}[row sep=0.6em]
0 \arrow[r] &
W_C^2 \arrow[r,"\Upsilon_2"] &
W_C^2 \arrow[r] &
\left.\widetilde{W_B^2/\Upsilon_2}\right|_{X_{\mathrm{fppf}}}
\arrow[r] &
0.
\end{tikzcd}
\end{equation*}
The associated long exact sequence in cohomology gives
\begin{equation*}
\begin{tikzcd}[row sep=0.6em]
0 \arrow[r] &
W_B^2(C) \arrow[r,"\Upsilon_2"] &
W_B^2(C) \arrow[r] &
\widetilde{W_B^2/\Upsilon_2}(C) \arrow[r] &
H^1(X_{\mathrm{fppf}},W_C^2) \arrow[r] &
\cdots.
\end{tikzcd}
\end{equation*}
By Theorem~\ref{st:orig-5-9} with $R=C$ and $r=2$, we have $H^1(X_{\mathrm{fppf}},W_C^2)=0$. Hence
\begin{align*}
\widetilde{W_B^2/\Upsilon_2}(C)
&\simeq
\Coker\left[
\Upsilon_2:W_B^2(C)\rightarrow W_C^2(C)
\right].
\end{align*}
Using the canonical identification $W_B^2(C)\simeq W_C^2(C)$, which is compatible with $\Upsilon_2$, we obtain
\begin{align*}
\widetilde{W_B^2/\Upsilon_2}(C)
\simeq
\Coker\left[\Upsilon_2:W_B^2(C)\rightarrow W_B^2(C)\right]
=(W_B^2/\Upsilon_2)(C).
\end{align*}
Since $C$ is arbitrary, the canonical morphism
\begin{align*}
W_B^2/\Upsilon_2
\rightarrow
\widetilde{W_B^2/\Upsilon_2}
\end{align*}
is an isomorphism. Therefore $W_B^2/\Upsilon_2$ is an fppf sheaf. The same argument applies to $W_B^2/\Upsilon_2'$.
\end{proof}
\subsection{Boundary cocycles for formal groups}\label{subsec:boundary-cocycles}

In this subsection, let $R$ be a commutative ring, and put
\begin{align*}
S=\Spec R,\qquad \mathcal H=\widehat{\GG}_{m,R},\qquad
\mathcal H(D)=1+\Nil(D)
\end{align*}
for an $R$-algebra $D$. We work in the category of abelian fppf
sheaves on $R$-algebras. Formal groups are regarded as sheaves of
nilpotent-valued points. All $\Ext^1$ groups in this subsection are
taken in this category.

For a commutative group sheaf $L$, put
\begin{align*}
C^1_0(L,\mathcal H)
&=\{a:L\to\mathcal H\mid a(e_L)=1\},\\
C^2_0(L,\mathcal H)
&=\{c:L\times_S L\to\mathcal H\mid c(e_L,x)=c(x,e_L)=1\},\\
Z^2_0(L,\mathcal H)
&=\left\{c\in C^2_0(L,\mathcal H)\ \middle|\
\begin{aligned}
c(x,y)&=c(y,x),\\
c(x,y)c(xy,z)&=c(y,z)c(x,yz)
\end{aligned}\right\}.
\end{align*}
The arrows in $C^1_0$ and $C^2_0$ mean morphisms of the underlying
sheaves. The group laws are pointwise multiplication, and the action
on $\mathcal H$ is trivial. Set
\begin{align*}
(\delta a)(x,y)&=\frac{a(x)a(y)}{a(xy)},\\
B^2_0(L,\mathcal H)&=\delta C^1_0(L,\mathcal H),\qquad
H^2_0(L,\mathcal H)=Z^2_0(L,\mathcal H)/B^2_0(L,\mathcal H).
\end{align*}
Then $\delta a\in Z^2_0(L,\mathcal H)$, since
\begin{align*}
(\delta a)(x,y)(\delta a)(xy,z)
=\frac{a(x)a(y)a(z)}{a(xyz)}
=(\delta a)(y,z)(\delta a)(x,yz).
\end{align*}
In particular, $\delta a=1$ if and only if $a\in\Hom(L,\mathcal H)$.

For $c\in Z^2_0(L,\mathcal H)$, let $P_c=\mathcal H\times_S L$
with product
\begin{align}\label{eq:boundary-extension-law}
(z,x)(z',y)=(zz'c(x,y)^{-1},xy).
\end{align}
This gives an extension of $L$ by $\mathcal H$. If $c'=c\delta a$, then
\begin{align*}
P_c\xrightarrow{\sim}P_{c'},\qquad(z,x)\longmapsto(za(x),x)
\end{align*}
is an isomorphism of extensions. Every such isomorphism has this form,
and the Baer sum of $P_c$ and $P_d$ is $P_{cd}$. Thus we have a natural
injective homomorphism
\begin{align*}
\jmath_L:H^2_0(L,\mathcal H)\hookrightarrow\Ext^1(L,\mathcal H),
\qquad[c]\longmapsto[P_c].
\end{align*}
Its image consists of the extensions whose underlying torsors have a
section. Indeed, a normalized section $s$ defines $c$ by
$s(x)s(y)=c(x,y)^{-1}s(xy)$, and $(z,x)\mapsto zs(x)$ identifies
$P_c$ with the given extension. For extensions and the Baer sum, see
\cite[Section~3.4]{Wei}.

We first prove a lemma on torsors. We use the following vanishing for
an $R_0$-module $M$:
\begin{align*}
H^1((\Spec R_0)_{\mathrm{fppf}},\widetilde M)=0.
\end{align*}
To see this, choose a finite affine fppf covering trivializing a
$\widetilde M$-torsor, and write its disjoint union as $\Spec D$.
The augmented Amitsur complex
\begin{align*}
0\longrightarrow M\longrightarrow M\otimes_{R_0}D
\longrightarrow M\otimes_{R_0}D^{\otimes2}
\longrightarrow M\otimes_{R_0}D^{\otimes3}\longrightarrow\cdots
\end{align*}
is exact. After tensoring with $D$, a contracting homotopy is given
by multiplying the extra $D$-factor with the first $D$-factor.
Faithful flatness gives exactness over $R_0$. Thus the descent
$1$-cocycle is a coboundary, and the local sections can be glued.

\begin{lem}\label{lem:pointed-formal-torsor}
Let $d\geq1$ and $X=\Spf R[[T_1,\ldots,T_d]]$, pointed by
$T_1=\cdots=T_d=0$. Let $P\to X$ be an $\mathcal H$-torsor in
the category of fppf sheaves. If a section $s_0$ over $S\to X$ is
given, then there exists a section $s:X\to P$ extending $s_0$.
\end{lem}
\begin{proof}
Put
\begin{align*}
I&=(T_1,\ldots,T_d),&
R_m&=R[[T_1,\ldots,T_d]]/I^{m+1},\\
X_m&=\Spec R_m,&
J_m&=I^m/I^{m+1}\subset R_m\quad(m\geq1).
\end{align*}
Then $J_m^2=0$ and $R_m/J_m=R_{m-1}$. Set $P_m=P\times_X X_m$.
Suppose that $s_{m-1}$ has been constructed. On the small fppf site
of $X_m$, consider the sheaf of lifts
\begin{align*}
\mathcal L_m(D)=
\{\sigma\in P_m(D)\mid
\sigma\bmod J_mD=s_{m-1}|_{\Spec(D/J_mD)}\},
\end{align*}
where $D$ is a flat $R_m$-algebra of finite presentation.

On a covering which trivializes $P_m$, the prescribed reduction is
$1+a_0\in\mathcal H(D/J_mD)$. For any lift $a\in D$ of $a_0$,
\begin{align*}
a_0^k=0\quad\Longrightarrow\quad a^k\in J_mD
\quad\Longrightarrow\quad a^{2k}=0.
\end{align*}
Hence $1+a\in\mathcal H(D)$ gives a local lift. Moreover,
\begin{align*}
\Ker[\mathcal H(D)\to\mathcal H(D/J_mD)]&=1+J_mD,\\
(1+b)(1+b')&=1+(b+b'),\\
J_m\otimes_{R_m}D&\xrightarrow{\sim}1+J_mD,
\qquad b\longmapsto1+b.
\end{align*}
The last isomorphism uses the flatness of $D$ over $R_m$.
Thus $\mathcal L_m$ is a $\widetilde J_m$-torsor, and
\begin{align*}
H^1((X_m)_{\mathrm{fppf}},\widetilde J_m)=0
\quad\Longrightarrow\quad
\mathcal L_m(R_m)\ne\varnothing.
\end{align*}
Choose $s_m\in\mathcal L_m(R_m)$. Induction gives
\begin{align*}
s_m|_{X_{m-1}}=s_{m-1}\qquad(m\geq1).
\end{align*}
For any $R$-algebra $D$, a point of $X(D)$ is a finite tuple of
nilpotent elements. These elements generate a nilpotent ideal, so
the point factors through some $X_m$. Therefore
\begin{align*}
X=\varinjlim_m X_m,\qquad s=(s_m)_m:X\longrightarrow P
\end{align*}
defines a section as a morphism of sheaves, with $s|_S=s_0$.
\end{proof}

\begin{thm}[Boundary-cocycle theorem]\label{thm:boundary-cocycle}
Let
\begin{align}\label{eq:boundary-general-sequence}
0\longrightarrow N\xrightarrow{\iota}G
\xrightarrow{\Psi}Q\longrightarrow0
\end{align}
be an exact sequence of commutative fppf sheaves over $R$.
Assume that $Q$ is a formal group whose underlying pointed formal
scheme is $\Spf R[[T_1,\ldots,T_d]]$ for some $d\geq1$.
Then
\begin{align*}
\jmath_Q:H^2_0(Q,\mathcal H)\xrightarrow{\sim}\Ext^1(Q,\mathcal H).
\end{align*}
For every $\chi\in\Hom(N,\mathcal H)$, there exist
\begin{align*}
c_\chi\in Z^2_0(Q,\mathcal H),\qquad
E_\chi\in C^1_0(G,\mathcal H)
\end{align*}
such that
\begin{align}\label{eq:boundary-general-cocycle}
\begin{split}
\partial_{\rm Ext}(\chi)&=\jmath_Q([c_\chi]),\\
\Psi^*c_\chi&=\delta E_\chi,\qquad \iota^*E_\chi=\chi,\\
E_\chi(g\iota(n))&=E_\chi(g)\chi(n).
\end{split}
\end{align}
The class $[c_\chi]$ is independent of the choices. The homomorphism
\begin{align*}
\partial_H:\Hom(N,\mathcal H)\longrightarrow H^2_0(Q,\mathcal H),
\qquad\chi\longmapsto[c_\chi]
\end{align*}
gives the exact sequence
\begin{equation}\label{eq:boundary-general-exact}
\begin{aligned}
1\longrightarrow\Hom(Q,\mathcal H)
&\xrightarrow{\Psi^*}\Hom(G,\mathcal H)
\xrightarrow{\iota^*}\Hom(N,\mathcal H)\\
&\xrightarrow{\partial_H}H^2_0(Q,\mathcal H)
\xrightarrow{\Psi^*}H^2_0(G,\mathcal H).
\end{aligned}
\end{equation}
\end{thm}
\begin{proof}
Let $0\to\mathcal H\to P\to Q\to0$ be an extension.
The identity of $P$ gives a section over the zero section of $Q$.
By Lemma~\ref{lem:pointed-formal-torsor}, the underlying torsor
has a normalized section. Hence $[P]\in\Im\jmath_Q$, proving
that $\jmath_Q$ is an isomorphism.

The boundary $\partial_{\rm Ext}(\chi)$ is represented by the
pushout of \eqref{eq:boundary-general-sequence} along $\chi$;
see \cite[Section~3.4]{Wei}. Explicitly, put
\begin{align*}
P_\chi&=(\mathcal H\times_SG)/N,&
(z,g)\cdot n&=(z\chi(n)^{-1},g\iota(n)),\\
j_\chi(z)&=[z,e_G],&
q_\chi([z,g])&=\Psi(g).
\end{align*}
The quotient is taken in fppf sheaves, and
\begin{align*}
0\longrightarrow\mathcal H\xrightarrow{j_\chi}P_\chi
\xrightarrow{q_\chi}Q\longrightarrow0
\end{align*}
is exact. In particular, $q_\chi$ is an $\mathcal H$-torsor.
Indeed,
\begin{align*}
\mathcal H\times_SG&\xrightarrow{\sim}P_\chi\times_QG,\\
(z,g)&\longmapsto([z,g],g).
\end{align*}
Choose a normalized section $s:Q\to P_\chi$ by
Lemma~\ref{lem:pointed-formal-torsor}. Define $c_\chi$ and $E_\chi$ by
\begin{align}\label{eq:boundary-section-formulas}
\begin{split}
s(x)s(y)&=j_\chi(c_\chi(x,y)^{-1})s(xy),\\
[1,g]&=j_\chi(E_\chi(g))s(\Psi(g)).
\end{split}
\end{align}
These are uniquely determined morphisms of sheaves. The group law
on $P_\chi$ gives
\begin{align*}
c_\chi(e,x)=c_\chi(x,e)&=1,\qquad c_\chi(x,y)=c_\chi(y,x),\\
c_\chi(x,y)c_\chi(xy,z)&=c_\chi(y,z)c_\chi(x,yz).
\end{align*}
Thus $c_\chi\in Z^2_0(Q,\mathcal H)$. The map
\begin{align*}
P_{c_\chi}\xrightarrow{\sim}P_\chi,
\qquad(z,x)\longmapsto j_\chi(z)s(x)
\end{align*}
is an isomorphism of extensions, so
$\partial_{\rm Ext}(\chi)=\jmath_Q([c_\chi])$.

Next, $E_\chi(e_G)=1$, and the quotient relation gives
\begin{align*}
[1,g\iota(n)]&=[\chi(n),g]
\quad\Longrightarrow\quad E_\chi(g\iota(n))=E_\chi(g)\chi(n),\\
E_\chi(\iota(n))&=\chi(n).
\end{align*}
Since $[1,g][1,h]=[1,gh]$, we also have
\begin{align*}
E_\chi(g)E_\chi(h)c_\chi(\Psi(g),\Psi(h))^{-1}=E_\chi(gh).
\end{align*}
This is $\Psi^*c_\chi=\delta E_\chi$.
If $s'=j_\chi(a)s$ for $a\in C^1_0(Q,\mathcal H)$, then
\begin{align*}
c'_\chi=c_\chi(\delta a)^{-1},\qquad
E'_\chi=E_\chi/(a\circ\Psi).
\end{align*}
Hence $[c'_\chi]=[c_\chi]$. Moreover,
$\partial_H=\jmath_Q^{-1}\partial_{\rm Ext}$ is a homomorphism.

The long exact $\Ext$ sequence gives exactness of
\eqref{eq:boundary-general-exact} at the three $\Hom$ terms.
Since $\jmath_G$ is injective and
$\jmath_G\Psi^*=\Psi^*\jmath_Q$, it also gives
\begin{align*}
\Ker[\Psi^*:H^2_0(Q,\mathcal H)\to H^2_0(G,\mathcal H)]
&=\jmath_Q^{-1}(\Im\partial_{\rm Ext})\\
&=\Im\partial_H.
\end{align*}
This proves the assertion.
\end{proof}

\begin{cor}\label{cor:boundary-comparison}
Under the assumptions of Theorem~\ref{thm:boundary-cocycle}, let
\begin{align*}
\chi\in\Hom(N,\mathcal H),\qquad
E\in C^1_0(G,\mathcal H),\qquad c\in Z^2_0(Q,\mathcal H).
\end{align*}
If $\iota^*E=\chi$ and $\delta E=\Psi^*c$, then
\begin{align*}
\partial_H(\chi)=[c].
\end{align*}
\end{cor}
\begin{proof}
Since $c(\Psi(g),e)=1$, we have
\begin{align*}
\frac{E(g)\chi(n)}{E(g\iota(n))}
=(\delta E)(g,\iota(n))=1.
\end{align*}
Thus the map
\begin{align*}
P_\chi\longrightarrow P_c,\qquad
[z,g]\longmapsto(zE(g),\Psi(g))
\end{align*}
is well-defined. The equality $\delta E=\Psi^*c$ shows that it is
a homomorphism. It induces the identity on $\mathcal H$ and $Q$,
so it is an isomorphism of extensions. Hence
$\partial_{\rm Ext}(\chi)=\jmath_Q([c])$, as desired.
\end{proof}

\begin{rem}\label{rem:boundary-generality}
Theorem~\ref{thm:boundary-cocycle} does not require $N$ to be finite
or $G$ to be representable. The integer $d$ is arbitrary, and no
Witt vector is chosen in its proof. The construction commutes with
base change: a normalized section and
\eqref{eq:boundary-section-formulas} pull back to the same formulas.
Thus the theorem can be used for the exact sequences in
\cite{A1,AA,Aki1}, or in higher dimension, whenever the sequence
is fppf exact and its quotient has the stated formal coordinates.
\end{rem}

\subsection{Proof of Theorem~\ref{st:orig-1-5}}
Let $B$ be an $A$-algebra such that each $\nu_{i,k}$ is nilpotent, and let $C$ be an arbitrary $B$-algebra. By Lemma~\ref{st:orig-4-15}, the following exact sequence is induced by the homomorphism $\Psi:=\Psi_{\uu,\vv}^{(l,j)}$:
\begin{equation}\label{eq:orig-5-6}
\begin{tikzcd}[row sep=0.6em]
0 \arrow[r] &
\mathcal{N}_C\arrow[r,"\iota"] &
\widehat{\EE}_C^{(\lambda_1,\lambda_2;D)} \arrow[r,"\Psi"] &
\widehat{\EE}_C^{(\lambda_1^{p^l},\lambda_2^{p^l};D')} \arrow[r] &
0,
\end{tikzcd}
\end{equation}
where $\iota$ is the canonical inclusion and the displayed $\Psi$ denotes the
formal homomorphism $\widehat{\Psi}_{\uu,\vv}^{(l,j)}$.  For brevity, put
\begin{align*}
G:=\widehat{\EE}_C^{(\lambda_1,\lambda_2;D)}
\quad\text{and}\quad
G':=\widehat{\EE}_C^{(\lambda_1^{p^l},\lambda_2^{p^l};D')}.
\end{align*}
The sequence \eqref{eq:orig-5-6} is exact as a sequence of fppf
sheaves on $C$-algebras. As pointed formal schemes, we have
\begin{align*}
G\simeq\Spf C[[X,Y]],\qquad G'\simeq\Spf C[[X,Y]].
\end{align*}
Apply Theorem~\ref{thm:boundary-cocycle} with
\begin{align*}
R=C,\qquad N=\mathcal N_C,\qquad Q=G',\qquad d=2.
\end{align*}
It gives a homomorphism
\begin{align*}
\partial_H:\Hom(\mathcal N_C,\widehat{\GG}_{m,C})
\longrightarrow H^2_0(G',\widehat{\GG}_{m,C}),\qquad
\jmath_{G'}\partial_H=\partial_{\rm Ext},
\end{align*}
and the exact sequence
\begin{equation}\label{eq:orig-5-7}
\begin{tikzcd}[row sep=0.6em, column sep=3.5ex]
1 \arrow[r] &
\Hom(G',\widehat{\GG}_{m,C}) \arrow[r, "(\Psi)^\ast"] &
\Hom(G,\widehat{\GG}_{m,C}) \arrow[r, "(\iota)^\ast"] &
\Hom(\mathcal{N}_C,\widehat{\GG}_{m,C})\\
  \arrow[r, "\partial_H"] &
H^2_0(G',\widehat{\GG}_{m,C}) \arrow[r, "(\Psi)^\ast"] &
H^2_0(G,\widehat{\GG}_{m,C}).
\end{tikzcd}
\end{equation}
On the other hand, by Lemma~\ref{st:orig-5-5} and Lemma~\ref{st:orig-5-7}, we have the commutative diagram
\begin{equation}\label{eq:orig-5-8}
\begin{tikzcd}
0 \arrow[r] &
W_B^2(C) \arrow[r,"\Upsilon_2"] \arrow[d,"U_2'"] &
W_B^2(C) \arrow[r,"\pi"] \arrow[d,"U_2"] &
(W_B^2/\Upsilon_2)(C) \arrow[r] \arrow[d,"\overline{U_2}"] &
0\\
0 \arrow[r] &
W_B^2(C) \arrow[r,"\Upsilon_2'"] &
W_B^2(C) \arrow[r,"\pi"] &
(W_B^2/\Upsilon_2')(C) \arrow[r] &
0.
\end{tikzcd}
\end{equation}
Here and below, $\Coker(U_2)(C)$ denotes the cokernel of $U_2(C)$ in the category of abelian groups, and similarly for $U_2'$.
The snake lemma applied to \eqref{eq:orig-5-8} gives the exact sequence
\begin{equation}\label{eq:orig-5-9}
\begin{tikzcd}[row sep=0.6em]
 &
\Ker(U_2')(C) \arrow[r,"\Upsilon_2"] &
\Ker(U_2)(C) \arrow[r,"\pi"] &
\Ker(\overline{U_2})(C)\\
\arrow[r,"\partial_W"] &
\Coker(U_2')(C) \arrow[r,"\Upsilon_2'"] &
\Coker(U_2)(C).
\end{tikzcd}
\end{equation}
Combining the exact sequences \eqref{eq:orig-5-7} and \eqref{eq:orig-5-9} with the isomorphisms \eqref{eq:orig-3-2} and \eqref{eq:orig-3-3}, we obtain the following diagram.  The third row continues the first, and the fourth continues the second.  All vertical maps except possibly $\phi_C$ are isomorphisms:
\begin{equation}\label{eq:orig-5-10}
\begin{tikzcd}[column sep=4.2ex]
\Hom(G',\widehat{\GG}_{m,C}) \arrow[r,"(\Psi)^\ast"] &
\Hom(G,\widehat{\GG}_{m,C}) \arrow[r,"(\iota)^\ast"] &
\Hom(\mathcal{N}_C,\widehat{\GG}_{m,C})\\
\Ker(U_2')(C) \arrow[r,"\Upsilon_2"] \arrow[u,"\phi_1"] & 
\Ker(U_2)(C) \arrow[r,"\pi"] \arrow[u,"\phi_2"] &
\mathcal{M}_B(C) \arrow[u,"\phi_C"]\\[-2ex]
\hspace*{10ex} \arrow[r,"\partial_H"] &
H^2_0(G',\widehat{\GG}_{m,C}) \arrow[r,"(\Psi)^\ast"] &
H^2_0(G,\widehat{\GG}_{m,C})\\
\hspace*{10ex} \arrow[r,"\partial_W"] &
\Coker(U_2')(C) \arrow[r,"\Upsilon_2'"] \arrow[u,"\phi_3"] &
\Coker(U_2)(C) \arrow[u,"\phi_4"],
\end{tikzcd}
\end{equation}
where we define $\phi_C$ using the correction cochains.
For $\overline{\vv}\in\mathcal M_B(C)$, choose a representative
$\vv\in W^2(C)$. There is a unique $\ww\in W^2(C)$ with
$U_2\vv=\Upsilon_2'\ww$: existence is the definition of
$\mathcal M_B(C)$, and uniqueness follows from
Lemma~\ref{st:orig-5-7}. Put
\begin{align}\label{eq:corrected-character}
\phi_C(\overline{\vv})
:=\iota^*\!\left(
E_p(\vv,(\lambda_1,\lambda_2);X,Y)\Gamma_{\ww}(X,Y)
\right).
\end{align}
We verify its character property and independence of the representative
below.

First establish the cocycle identity needed in this verification.
For every $\ww\in W^2(C)$,
\begin{align}\label{eq:orig-5-11}
\begin{split}
&F_p(\Upsilon_2'\ww,(\lambda_1,\lambda_2);\XX,\YY)\\
&\quad=
F_p(\ww,(\lambda_1^{p^l},\lambda_2^{p^l});\Psi(\XX),\Psi(\YY))
\frac{\Gamma_{\ww}(\XX\cdot\YY)}
{\Gamma_{\ww}(\XX)\Gamma_{\ww}(\YY)}.
\end{split}
\end{align}
Here $\Gamma_{\ww}(\XX)$ means evaluation of the same two-variable
series at $\XX=(X_0,X_1)$, and similarly for $\YY$.
To prove the identity, use Sublemma~\ref{st:orig-5-6} to choose a
faithfully flat $C$-algebra $C'$ and $\zz\in W^2(C')$ satisfying
$U_2'\zz=\ww|_{C'}$. Over $C'$, apply the coboundary to
\eqref{eq:corrected-pullback} and use \eqref{eq:orig-3-1} and
\eqref{eq:upsilon-factorization}. This gives
\begin{align*}
\Psi^*F_p(U_2'\zz,(\lambda_1^{p^l},\lambda_2^{p^l}))
&=F_p(U_2\Upsilon_2\zz,(\lambda_1,\lambda_2))
\,\delta\Gamma_{U_2'\zz}\\
&=F_p(\Upsilon_2'\ww,(\lambda_1,\lambda_2))
\,\delta\Gamma_{\ww}.
\end{align*}
Rearranging proves \eqref{eq:orig-5-11} over $C'$.
Both sides are formal series over $C$, and $C\to C'$ is injective,
so the identity descends coefficientwise to $C$.

\begin{rem}\label{st:orig-5-11}
We now verify that \eqref{eq:corrected-character} is a well-defined
homomorphism. For $U_2\vv=\Upsilon_2'\ww$, set
\begin{align*}
H_{\vv,\ww}(X,Y)
:=E_p(\vv,(\lambda_1,\lambda_2);X,Y)\Gamma_{\ww}(X,Y).
\end{align*}
By \eqref{eq:orig-3-1} and \eqref{eq:orig-5-11},
\begin{align}\label{eq:corrected-defect}
\delta H_{\vv,\ww}
=\Psi^*F_p(\ww,(\lambda_1^{p^l},\lambda_2^{p^l})).
\end{align}
The target cocycle is normalized. Restriction to
$\mathcal N_C\times_C\mathcal N_C$, using $\Psi\iota=0$, gives
$(\iota\times\iota)^*\delta H_{\vv,\ww}=1$.
Thus $\iota^*H_{\vv,\ww}$ is a character.

If the representative changes to $\vv+\Upsilon_2\zz$, its unique
witness changes to $\ww+U_2'\zz$. Additivity of $E_p$ and $\Gamma$,
together with \eqref{eq:corrected-pullback}, gives
\begin{align*}
H_{\vv+\Upsilon_2\zz,\,\ww+U_2'\zz}
&=H_{\vv,\ww}\,
E_p(\Upsilon_2\zz,(\lambda_1,\lambda_2))
\Gamma_{U_2'\zz}\\
&=H_{\vv,\ww}\,\Psi^*E_p(\zz,(\lambda_1^{p^l},\lambda_2^{p^l})).
\end{align*}
The final factor restricts to $1$ on $\mathcal N_C$.
This proves independence of the representative.
If two classes have representatives $\vv,\vv'$ with witnesses
$\ww,\ww'$, their sum has witness $\ww+\ww'$.
Additivity therefore also proves that $\phi_C$ is a homomorphism.
All ingredients commute with base change, and uniqueness makes the
witness commute with base change as well. Hence the maps $\phi_C$
form a morphism of fppf sheaves.
\end{rem}

\begin{lem}\label{st:orig-5-12}
Diagram \eqref{eq:orig-5-10}, with $\phi_C$ defined by
\eqref{eq:corrected-character}, is commutative. Consequently
$\phi_C$ is an isomorphism.
\end{lem}
\begin{proof}
The equality $\Psi^*\phi_1=\phi_2\Upsilon_2$ is
Lemma~\ref{st:orig-5-4}. If $\vv\in\Ker U_2$, its witness in
\eqref{eq:corrected-character} is $0$, since $\Upsilon_2'$ is
injective. As $\Gamma_0=1$, this also gives
$\iota^*\phi_2=\phi_C\pi$.

The last factor of \eqref{eq:orig-5-11} is the inverse of the
coboundary of the normalized cochain $\Gamma_{\ww}$ on $G$.
Consequently,
\begin{align*}
[F_p(\Upsilon_2'\ww,(\lambda_1,\lambda_2))]
=[\Psi^*F_p(\ww,(\lambda_1^{p^l},\lambda_2^{p^l}))]
\quad\text{in }H^2_0(G,\widehat{\GG}_{m,C}).
\end{align*}
This proves $\Psi^*\phi_3=\phi_4\Upsilon_2'$.

It remains to compare the two boundary maps.
Take $\overline{\vv}\in\mathcal M_B(C)$, with representative $\vv$
and witness $\ww$. The snake-lemma boundary is
\begin{align*}
\partial_W(\overline{\vv})=[\ww]\in\Coker(U_2')(C).
\end{align*}
Put
\begin{align*}
H&:=H_{\vv,\ww}\in C^1_0(G,\widehat{\GG}_{m,C}),&
\chi&:=\iota^*H=\phi_C(\overline{\vv}),\\
C_{\ww}&:=F_p(\ww,(\lambda_1^{p^l},\lambda_2^{p^l}))
   \in Z^2_0(G',\widehat{\GG}_{m,C}).
\end{align*}
By \eqref{eq:corrected-defect}, we have
\begin{align*}
\iota^*H=\chi,\qquad \delta H=\Psi^*C_{\ww}.
\end{align*}
Apply Corollary~\ref{cor:boundary-comparison} with $E=H$ and
$c=C_{\ww}$. We obtain
\begin{align*}
\partial_H(\phi_C(\overline{\vv}))
=[C_{\ww}]
=[F_p(\ww,(\lambda_1^{p^l},\lambda_2^{p^l}))]
=\phi_3(\partial_W(\overline{\vv})).
\end{align*}
All squares of \eqref{eq:orig-5-10} now commute.
The two five-term sequences are exact, and the four vertical maps other than $\phi_C$ are
isomorphisms by \eqref{eq:orig-3-2} and \eqref{eq:orig-3-3}.
The five lemma proves that $\phi_C$ is an isomorphism.
\end{proof}

Finally, by \eqref{eq:orig-4-8},
\begin{align*}
\mathcal M_B(C)
&\simeq\Hom(\mathcal N_C,\widehat{\GG}_{m,C})\\
&=\Hom(\mathcal N_C,\GG_{m,C})
=(\mathcal N_B)^D(C).
\end{align*}
The isomorphisms are natural in $C$, proving
Theorem~\ref{st:orig-1-5}.

\subsection{Some corollaries}
We first consider the case of positive characteristic.
\begin{cor}[{cf. \cite[Theorem~2]{AA}}]\label{st:orig-5-13}
Assume that $A$ is an $\FF_p$-algebra. Let $\lambda_1$ and $\lambda_2$ be nonzero divisors of $A$, and assume that the image of $\lambda_1$ in $A_{\lambda_2}$ is nilpotent. Let $j=0$, and choose $\uu$, $\vv$, $D$, and $D'$ as in Remark~\ref{st:orig-4-13}.
Choose a lift $\uu^\#\in W(A)$ of the chosen $\widetilde{\uu}$ and set
\begin{align*}
\widetilde{\vv}:=F^l(\uu^\#)
\end{align*}
when defining the Witt-vector operators, as in Remark~\ref{st:orig-5-8}.
Then
\begin{align*}
\Psi_{\uu,\vv}^{(l,0)}(X,Y)=\left(X^{p^l},Y^{p^l}\right)\quad\text{and}\quad
\Upsilon_2=\Upsilon_2'=V^l\times V^l.
\end{align*}
Consequently, by {\rm \eqref{eq:orig-1-10}}, we have an isomorphism
\begin{align*}
\left(\mathcal{N}_A\right)^D
\simeq\mathcal{M}_A=\Ker\!\left[\overline{U_2}:W_{l,A}^2\longrightarrow W_{l,A}^2\right].
\end{align*}
\end{cor}
\begin{proof}
Since $A$ is an $\FF_p$-algebra, we may take $\nu_{i,k}=0$ for each $k$. Hence every $\nu_{i,k}$ is nilpotent in $A$, and therefore Theorem~\ref{st:orig-1-5} applies. By Remark~\ref{st:orig-4-13}, with the above choices of $\vv$, $D$, and $D'$, we have
\begin{align*}
\Psi_{\uu,\vv}^{(l,0)}(X,Y)=\left(X^{p^l},Y^{p^l}\right).
\end{align*}
Moreover, by Remark~\ref{st:orig-5-8}, we obtain
\begin{align*}
\Upsilon_2=\Upsilon_2'=V^l\times V^l.
\end{align*}
Since $W_A/V^l\simeq W_{l,A}$, the homomorphism induced by $U_2$ on the quotient $W_A^2/(V^l\times V^l)$ is identified with $\overline{U_2}:W_{l,A}^2\rightarrow W_{l,A}^2$.
\end{proof}
\begin{rem}\label{st:orig-5-14}
{\rm
The same conclusion holds after any base change $A\to B$. Indeed, since $A$
is an $\FF_p$-algebra, all the elements $\nu_{i,k}$ may be chosen to be zero,
and their images in $B$ are therefore nilpotent. Thus Theorem~\ref{st:orig-1-5} applies to
the base-changed isogeny. In particular, although $\lambda_1$ and $\lambda_2$
are assumed to be nonzero divisors in $A$, their images in $B$ may be
nilpotent. For example, for positive integers $N_1$ and $N_2$, one may take
\begin{align*}
B=A/(\lambda_1^{N_1},\lambda_2^{N_2}).
\end{align*}
For this base change, we obtain
\begin{align*}
\left(\mathcal N_B\right)^D
\simeq
\mathcal M_B
=\Ker\!\left[\overline{U_2}:W_{l,B}^2\longrightarrow W_{l,B}^2\right].
\end{align*}
Here all the objects over $B$ are understood as the base changes of the
corresponding objects over $A$.
}
\end{rem}
We next apply Theorem~\ref{st:orig-1-5} to quotients of cyclotomic integer rings.
\begin{cor}\label{st:orig-5-15}
Let $m$, $m_1$, and $m_2$ be positive integers such that $l<m_i\leq m$ $(i=1,2)$. Let $\zeta_{p^m}$ be a primitive $p^m$-th root of unity, and put
\begin{align*}
A:=\ZZ_{(p)}[\zeta_{p^m}],\quad
\lambda:=\zeta_{p^m}-1\quad \text{and}\quad
\lambda_i:=\zeta_{p^{m_i}}-1
\end{align*}
for $i=1,2$. For a positive integer $n$, put $\mu=\lambda^n$ and set $A_\mu:=A/(\mu)$. Assume that $\uu$, $\vv$, $D$, and $D'$ satisfy the assumptions of Theorem~\ref{st:orig-1-2} over $A$. Then we have
\begin{align*}
\mathcal{M}_{A_\mu}
\simeq
\left(\mathcal{N}_{A_\mu}\right)^D.
\end{align*}
\end{cor}
\begin{proof}
Let $v$ be the valuation of $A$, normalized by $v(\lambda)=1$. Then $v(p)=p^m-p^{m-1}=p^{m-1}(p-1)$ and $v(\lambda_i)=p^{m-m_i}$ $(i=1,2)$ hold. For $0\leq k\leq l-1$, write $p^{l-k}=\nu_{i,k}\lambda_i^{p^l-p^k}$. Then the equalities
\begin{align*}
v(\nu_{i,k})
=(l-k)v(p)-(p^l-p^k)v(\lambda_i)
=p^{m-m_i}\left\{(l-k)p^{m_i-1}(p-1)-(p^l-p^k)\right\}
\end{align*}
hold. Since $p^l-p^k=(p-1)\sum_{r=k}^{l-1}p^r$ and $m_i>l$, we have
\begin{align*}
(l-k)p^{m_i-1}(p-1)
>
(p-1)\sum_{r=k}^{l-1}p^r
=
p^l-p^k.
\end{align*}
Hence, since $v(\nu_{i,k})>0$, we obtain $\nu_{i,k}\in(\lambda)$ for every $i=1,2$ and $0\leq k\leq l-1$. Since the image of $\lambda$ in $A_\mu$ is nilpotent, the image of every $\nu_{i,k}$ in $A_\mu$ is nilpotent. Therefore the
assumptions of Theorem~\ref{st:orig-1-5} are satisfied, and the assertion follows.
\end{proof}
\begin{rem}\label{st:orig-5-16}
{\rm
The strict inequalities $l<m_1$ and $l<m_2$ in Corollary~\ref{st:orig-5-15} are essential for applying Theorem~\ref{st:orig-1-5}. Indeed, if $m_i=l$, then $\nu_{i,l-1}$ is a unit in $A$, as observed in Remark~\ref{st:orig-4-1}, and hence its image in the quotient $A_{\mu}$ remains a unit and is not nilpotent.
}
\end{rem}

\section{Declaration of generative AI and AI-assisted technologies in the manuscript preparation process}

During the preparation of this work, the author used GPT-6 Astra to assist with reviewing and revising the manuscript. The author reviewed and edited the output as needed and takes full responsibility for the content of the published article. 

\section*{Acknowledgements}
The author would like to express his sincere gratitude to Professor~Tsutomu~Sekiguchi for his valuable suggestions. He is deeply grateful to Professor~Akira~Masuoka for his generous support and warm encouragement. He also thanks Professor~Takahiro~Noi for his helpful advice on improving the presentation of this paper.

\bibliographystyle{amsplain}
\bibliography{ref_amano}

\providecommand{\bysame}{\leavevmode\hbox to3em{\hrulefill}\thinspace}
\providecommand{\MR}{\relax\ifhmode\unskip\space\fi MR }
\providecommand{\MRhref}[2]{%
  \href{http://www.ams.org/mathscinet-getitem?mr=#1}{#2}
}
\providecommand{\href}[2]{#2}
\begin{thebibliography}{10}

\bibitem{Aki1}
N.~Aki, \emph{\textrm{On the Cartier duality of certain finite group schemes of
  type $(p^n,\ldots,p^n)$}}, Tsukuba J. Math. \textbf{34} (2010), no.~2,
  183--200.

\bibitem{AA}
N.~Aki and M.~Amano, \emph{\textrm{On the Cartier duality of certain finite
  group schemes of type $(p^n,p^n)$}}, Tsukuba J. Math. \textbf{34} (2010),
  no.~1, 31--46.

\bibitem{A1}
M.~Amano, \emph{\textrm{On the Cartier duality of certain finite group schemes
  of order $p^n$}}, Tokyo J. Math. \textbf{33} (2010), 117--127.

\bibitem{A4}
\bysame, \emph{\textrm{On the injectivity of certain homomorphisms between
  extensions of $\widehat{\mathcal{G}}^{(\lambda)}$ by $\widehat{\mathbb{G}}_m$
  over a $\mathbb{Z}_{(p)}$-algebra}}, Tsukuba J. Math. \textbf{49} (2025),
  no.~2, 255--275.

\bibitem{DG}
M.~Demazure and P.~Gabriel, \emph{Groupes alg\'{e}briques}, Tome I,\
  Masson-North-Holland, Paris-Amsterdam, 1970.

\bibitem{GM}
B.~Green and M.~Matignon, \emph{\textrm{Liftings of Galois Covers of Smooth
  Curves}}, Compositio Mathematica \textbf{113} (1998), no.~3, 237--272.

\bibitem{HZ}
M.~Hazewinkel, \emph{Formal~groups~and~applications}, Academic~Press, New~York,
  1978.

\bibitem{Jan}
Uwe Jannsen, \emph{Continuous {\'e}tale cohomology}, Mathematische Annalen
  \textbf{280} (1988), no.~2, 207--245.

\bibitem{MRT1}
A.~M\'{e}zard, M.~Romagny, and D.~Tossici, \emph{\textrm{Models of group
  schemes of roots of unity}}, Ann. Inst. Fourier, Grenoble \textbf{63} (2013),
  no.~3, 1055--1135.

\bibitem{MRT2}
\bysame, \emph{\textrm{Sekiguchi-Suwa theory revisited}}, Journal de
  Th\'{e}orie des Nombres de Bordeaux \textbf{26} (2014), 163--200.

\bibitem{S1}
T.~Sekiguchi, \emph{\textrm{On the deformations of Witt groups to tori, II}},
  Journal of Algebra \textbf{138} (1991), 273--297.

\bibitem{SOS}
T.~Sekiguchi, F.~Oort, and N.~Suwa, \emph{\textrm{On the deformation of
  Artin-Schreier to Kummer}}, Ann. Sci. \'{E}cole Norm. Sup. (4) \textbf{22}
  (1989), no.~3, 345--375.

\bibitem{SS2}
T.~Sekiguchi and N.~Suwa, \emph{\textrm{On the unified
  Kummer-Artin-Schreier-Witt theory}}, Pr\'{e}publication No.~111,
  Universit\'{e} de Bordeaux (1999).

\bibitem{SS1}
\bysame, \emph{\textrm{A note on extensions of algebraic and formal groups
  IV}}, Tohoku Math. J \textbf{53} (2001), 203--240.

\bibitem{SL}
J.~P. Serre, \emph{Local fields}, Springer, New~York, 1979.

\bibitem{SAG}
\bysame, \emph{Algebraic groups and class fields}, Springer, New~York, 1988.

\bibitem{To}
D.~Tossici, \emph{\textrm{Models of $\mu_{p^2,K}$ over discrete valuation
  ring}}, Journal of Algebra \textbf{323} (2010), 1908--1957.

\bibitem{W}
W.~Waterhouse, \emph{\textrm{A Unified Kummer-Artin-Schreier Sequence}},
  Mathematische Annalen \textbf{277} (1987), 447--452.

\bibitem{WW}
W.~Waterhouse and B.~Weisfeiler, \emph{\textrm{One-dimensional affine group
  schemes}}, Journal of Algebra \textbf{66} (1980), 550--568.

\bibitem{Wei}
C.~Weibel, \emph{An introduction to homological algebra}, Cambridge University
  Press, Cambridge, 1994.

\end{thebibliography}
\end{document}